\documentclass[12pt]{article}
\usepackage{a4wide}

\usepackage{fancyhdr, amsmath, amssymb, amsfonts, url, mathrsfs, graphicx, epsfig,  amsthm}
\usepackage{tikz}

\newtheorem{thm}{Theorem}[section]
\newtheorem{cor}[thm]{Corollary}
\newtheorem{lemma}[thm]{Lemma}
\newtheorem{proposition}[thm]{Proposition}
\newtheorem{definition}[thm]{Definition}
\newtheorem{example}[thm]{Example}
\newtheorem{remark}[thm]{Remark}

\newcommand{\marginnotem}[1]{\mbox{}\marginpar{\tiny\raggedright\hspace{0pt}{\color{magenta} {\bf Mark}$\blacktriangleright$  #1}}}

\newcommand\eqdef{\buildrel{\rm def}\over=}
\newcommand\tr{\mathop{\rm tr}\nolimits }
\newcommand\dif{\mathord{\rm d}}
\newcommand\acosh{\mathop{\rm ArcCosh}\nolimits }

\usepackage{pgfplots}

\AtEndDocument{\bigskip{\footnotesize
  \noindent
  \textsc{D. Parmenter, Department of Mathematics, Warwick University, Coventry, CV4 7AL, UK}
  
  \noindent
  \textit{E-mail address}: \texttt{d.parmenter@warwick.ac.uk}
  \par
  \addvspace{\medskipamount}
  \noindent
  \textsc{M. Pollicott, Department of Mathematics, Warwick University, Coventry, CV4 7AL, UK.}
  
  \noindent
  \textit{E-mail address}: \texttt{masdbl@warwick.ac.uk}
}}

\begin{document}
 \title{Gibbs measures for  hyperbolic attractors \\ defined  by  densities}
\author{David Parmenter and Mark Pollicott
\thanks{
The second author is partly
supported by ERC-Advanced Grant 833802-Resonances and EPSRC grant
EP/T001674/1}
}
\maketitle

\abstract{In this article we will  describe a new construction for Gibbs measures for hyperbolic attractors generalizing the original construction of Sinai, Bowen and Ruelle of SRB measures.  The classical  construction
of the SRB measure is based on pushing forward the normalized volume on a piece of unstable manifold.   By modifying the density at each step appropriately we show that the resulting measure is a prescribed Gibbs measure.  This contrasts with, and complements,   the construction of Climenhaga-Pesin-Zelerowicz who replace the volume on the unstable manifold by a fixed reference measure.  Moreover,  the simplicity of our proof, which uses  only explicit properties on the growth rate of unstable manifold and entropy estimates,  has the additional advantage that it applies in more general settings.}

\section{Introduction}
In broad terms, smooth ergodic theory describes the study of invariant measures for diffeomorphisms of compact manifolds.   
  Gibbs measures form an especially  natural  family of  invariant probability measures  
  which have played   an important role, particularly in the study of hyperbolic dynamical systems,  for over   50 years.   The best known examples of such measures  include the  well known Sinai-Ruelle-Bowen measures and the measure of maximal entropy, also known as the Bowen-Margulis measure.
  More generally, Gibbs measures  are invariant probability measures $\mu_{G}$ associated to H\"older continuous potential functions 
  $G$.  In this article we will describe a simple construction based on the push forward of the induced  volume on local unstable manifolds.   This generalizes the original method  used  by Sinai and Ruelle to construct the SRB measure.

  \subsection{Gibbs measures for hyperbolic diffeomorpisms}
Given H\"older continuous function $G$  there are various ways to describe  the associated Gibbs measures.   One standard  approach due to Ruelle and Walters is to use the variational principle, which characterizes the Gibbs measures 
as {\it equilibrium states}.  More generally, given a continuous map  $f: X \to X$ and a continuous function $G: \textcolor{black}{X} \to \mathbb R$
we say that a\textcolor{black}{n} $f$-invariant probability measure  $\mu_G$ is an equilibrium state for $G$ if 
$$h(\mu_G) + \int G d\mu_G = \sup\left\{h(\mu) + \int G d\mu \hbox{ : } \mu \in \mathcal M_f(X) \right\} \eqno(1.1)$$
where $\mathcal M_f(X)$ is the space of $f$-invariant probability measures
i.e., $\mu_G$ is a measure which maximizes the sum of the entropy and the integral of $G$,  over all $f$-invariant measures \cite{ruelle}, \cite{walters}.
The value attained by  the supremum in (1.1)  is called the pressure and denoted $P(G)$.


It is known that for every continuous function $G: X \to \mathbb R$ 
{\color{black} and expansive map $f$}  there exists at least one equilibrium state.  Under the additional assumptions that $G$ is H\"older continuous 
{\color{black} and $f$ is uniformly hyperbolic}
then there is a unique measure $\mu_G \in \mathcal M_f(X)$ satisfying (1.\textcolor{black}{1}) which is the Gibbs measure for $G$.
 \cite{bowen-SLN}.
One   approach  to constructing  Gibbs measures is to use weighted Dirac measures on  periodic orbits, which was developed by Bowen \cite{bowen-periodic} and Parry \cite{parry1}, \cite{parry2}.
We refer the reader to the recent comprehensive survey of Climenhaga et al for a detailed  discussion of constructing Gibbs measures, including their own original    construction \cite{CPZ-1}.
{\color{black} For other perspectives on the construction of Gibbs measures we refer the reader to 
\cite{gibbs1}, \cite{gibbs2}, \cite{gibbs3} and \cite{gibbs4}. }


To explain the viewpoint we  take in this article we begin by recalling there  is a well known construction of the classical 
Sinai-Ruelle-Bowen measure 
$\mu_{SRB}$ using the push forward of the 
normalized 
volume $\lambda$ on any  piece of local unstable manifold $W^u_\delta(x)$. 
Of course it follows from the classical 
Krylov–Bogolyubov
theorem on the existence of invariant measures 
 that the weak star limit points of
$$\frac{1}{n}\sum_{k=0}^{n-1} f_*^k \lambda, \quad n \geq 1,$$
are $f$-invariant.
(Here we denote by  $f_*^k\lambda (A) = \lambda(f^{-k}A)$  the push forward measure supported on $f^kW^u_\delta(x)$.)
However, in the context of hyperbolic attractors  much more is true.
More precisely, there is the following  famous result due to Sinai (for the particular case of Anosov diffeomorphisms)
and Ruelle (in the general setting of hyperbolic attractors).



\begin{thm}[Sinai \cite{sinai}, Ruelle   \cite{ruelle}]\label{thm:sinai}
Let $f: X \to X$ be a $C^{1+\alpha}$ topologically mixing   hyperbolic attractor. 
Given $x\in X$ and $\delta > 0$ 
consider 
a (normalized) volume measure $\lambda = \lambda_{W_\delta^u(x)}$ on a piece of local unstable manifold $W_\delta^u(x)$.
Then 
 the averages 
 $$\mu_n^{SRB} = \frac{1}{n}\sum_{k=0}^{n-1} f_*^k \lambda, \quad n \geq 1,$$
  converge in the weak star topology to $\mu_{SRB}$ as $n \to +\infty$\footnote{Sinai and Ruelle actually show the stronger result that the  measures $f_*^k\lambda$ converge to $\mu_{SRB}$ in the weak star topology without the need to average.  Moreover, the topological mixing hypothesis is not restrictive because of the Smale spectral decomposition theorem \cite{smale}}.
\end{thm}


  In the case of the SRB measures, it seems  particularly  appropriate  to use $\lambda$ as a starting measure since the limiting measure  $\mu_{SRB}$ is an $f$-invariant measure which is supported on $X$, but induces an absolutely continuous probability measure on all of the unstable manifolds.  

The SRB measure can also be viewed as a 
particular example of a Gibbs measure associated to the  H\"older continuous function $\Phi: X \to \mathbb R$ defined by  
$$\Phi(x) = -  \log \left|\det (Df|E_x^u)\right|,$$
which is usually called the unstable expansion coefficient, where $E^u_x$ is the unstable bundle for the diffeomorphism. 
  More generally, given any  H\"older continuous function $G : X \to \mathbb R$  we present an  analogous construction of the  Gibbs measure $\mu_G$. 

\begin{thm}\label{thm:gibbsdiffeo}
Let $f: X \to X$ be a $C^{1+\alpha}$ topologically mixing  hyperbolic attracting diffeomorphism and  let  $G: X \to \mathbb R$
be a H\"older continuous function. 
Given $x\in X$ and $\delta > 0$ 
consider 
the  sequence of probability measures $(\lambda_n)_{n=1}^\infty$ 
supported on $W^u_\delta(x)$ 
and absolutely continuous with respect to the induced volume $\lambda = \lambda_{W^u_\delta(x)}$
with densities
$$
\frac{d\lambda_n}{d\lambda}(y) :=  
\frac{
 \exp\left( \sum_{i=0}^{n-1} (G - \Phi)(f^iy)) \right) 
 }{
 \int_{W_\delta^u(x)}
 \exp\left( \sum_{i=0}^{n-1} (G-\Phi)(f^iz) \right) 
 d\lambda(z)
 }
\quad  \hbox{ for } y \in  W^u_\delta(x).\eqno(1.2)
$$
Then the averages  
$$\mu_n :=\frac{1}{n}\sum_{k=0}^{n-1} f_*^k \lambda_n, \quad n \geq 1,
\eqno(1.3)
$$
 converge in the weak star topology to $\mu_G$.
\end{thm}

It is clear from the statement that the  construction is independent of the choice of $x$ and $\delta > 0$.   
{\color{black} 

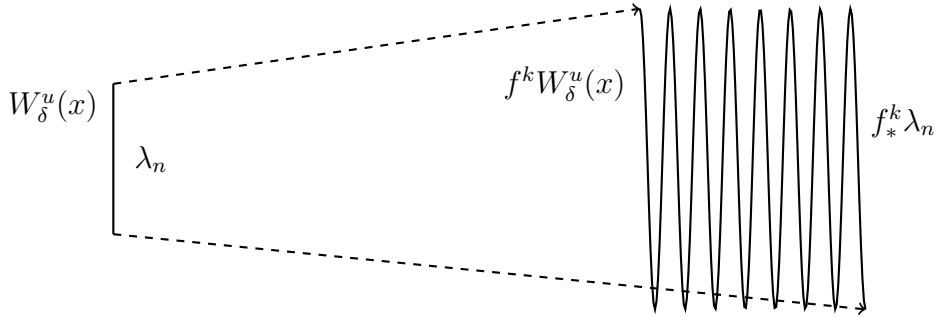
\begin{figure}[h!]
\centerline{
\begin{tikzpicture}[thick,scale=1.0, every node/.style={scale=1.0}]
\draw (-5,1) --  (-5,3);
\node at (-5.8,2.7) {$W_\delta^u(x)$};
\node at (1.0,3) {$f^{\textcolor{black}{k}}W_\delta^u(x)$};
\node at (-4.5,2.0) {$\lambda_{\textcolor{black}{n}}$};
\node at (5.5,2.5) {$f^{\textcolor{black}{k}}_*\lambda_{\textcolor{black}{n}}$};
\draw[dashed,->] (-5,1) -- (5,0);
\draw[dashed,->] (-5,3) -- (2,4);
\draw [domain= 2:5, samples=150] plot (\x, {2+2*cos(5*pi*\x r)});
\end{tikzpicture}
}
\caption{The  push forward of  the  measure $\lambda_n$ on $W_\delta^u(x)$ by $f^{\textcolor{black}{k}}$.}
\end{figure}

The proof  of Theorem \ref{thm:gibbsdiffeo}  actually gives  the following result.  

\begin{thm}\label{thm:gibbsdiffeo-1}
If $G: X \to \mathbb R$ is continuous then the weak star accumulation points for the measures (1.3) are 
equilibrium states for $G$.
\end{thm}
}

The term $\Phi$ appearing in (1.2) comes from the change of variables for
 $f^n: W^u_\delta(x) \to f^n(W^u_\delta(x))$ for which the Jacobian is 
$$\left|\det (Df^n|E_y^u)\right|
= \prod_{i=1}^n \left|\det  (Df|E_{f^iy}^u)\right| = 
 \exp\left( -\sum_{i=0}^{n-1}  \Phi(f^iy) \right)
 $$
 by the chain rule.  Thus we can reformulate (1.2) as 
 $$
 \lambda_n(A) = 
 \frac{\int_{f^{\textcolor{black}{n}}(A)} \exp \left( \sum_{i=1}^{n} G(f^{-i}y)\right) d\lambda_{f^kW_\delta^u}(y)}
 {\int_{f^nW_\delta^u(x)} \exp \left( \sum_{i=1}^{n} G(f^{-i}z)\right) d\lambda_{f^kW_\delta^u}(z)}
 \hbox{  for Borel  $A \subset W_\delta^u(x)$}
 $$ 
 which is   more convenient  in the proofs.

The proof of Theorem \ref{thm:gibbsdiffeo} will be given in \S 5. 
The simple nature of the proof suggests a number of potential generalizations, which we describe in \S 7.

An interesting aspect  of this approach  is that   all  of the 
reference measures $\lambda_n$ on $W^u_\delta(x)$  are absolutely continuous 
with respect to the volume $\lambda$,  although they now typically depend on $n$
(except, of course,  in the original  case of SRB measures  in  Theorem \ref{thm:sinai}
where the same measure $\lambda$ is used throughout).     This is reminiscent of the weighted standard pair construction  
in (5.12)  of \cite{jl} which was used to prove stronger statistical properties in a different setting.
In contrast to this, in a  recent very  interesting article, Climenhaga-Pesin-Zelerowicz use a different method to construct  Gibbs measures by pushing forward a fixed reference measure
on $W^u_\delta(x)$
 (dependent on $G$ and typically singular with respect to the leaf volume) constructed via a Carath\'eodory type construction \cite{CPZ-1}.

 \begin{remark}
 It is clear from the proof that the  conclusion of Theorem \ref{thm:gibbsdiffeo} remains true if $W^u_\delta(x)$ is replaced by an embedded disk $\mathcal D \subset M$ of dimension $E^u$
 provided it 
 {\color{black} is transverse to the strong stable manifolds.}
 \end{remark}
 
 \begin{remark}
The examples of    attractors  from \cite{williams} 
can be seen to  have an exponential rate  of convergence of $\int G d\mu_n^{SRB}$ to $\int G d\mu_{SRB}$ (for $G$ H\"older continuous)
 in Theorem 1.1, but this rate  can be made  arbitrarily small.
 \end{remark}
 
 {\color{black}
  \begin{remark}
  One of the advantages of considering attractors is that one can naturally study the induced volume on the unstable manifolds as a reference measure.  There is a prospect that these results can be extended to basic sets (see \cite{kh}, \S 6.4.d)  using another ambient measure,
  which will be pursued elsewhere.
  \end{remark}
  }

\subsection{Gibbs measures for hyperbolic  flows}
The results for diffeomorphisms in the previous subsection have analogues for hyperbolic attracting flows $\phi_t: X \to X$.
Given a H\"older continuous function $G: X \to \mathbb R$ we can again use a   variational principle,
to first  characterize the Gibbs measure $\mu_G$
as an equilibrium state. More generally, 
given a continuous flow $\phi_t:X \to X$
and a continuous function $G: \textcolor{black}{X} \to \mathbb R$ 
we say a $\phi$-invariant probability measure $\mu_G$ is an equilibrium state for $G$ if 
$$h(\mu_G) + \int G d\mu_G
 = \sup
 \left\{
 h(\mu) + \int G d\mu
  \hbox{ : } 
\mu \in \mathcal M_\phi(X)
\right\} 
\eqno(1.4)$$
where $\mathcal M_\phi(X)$ is the space of $\phi$-invariant probability measures and 
$h(\mu) $ is the entropy of the time-one flow $\phi_{t=1}: X \to X$ and $\mu \in \mathcal M_\phi\textcolor{black}{(X)}$, 
i.e., $\mu_G$ is a measure which maximizes the sum of the entropy and the integral over all $\phi$-invariant measures \cite{ruelle}, \cite{walters}.
Furthermore, the value of the supremum is again called the pressure
and denoted  $P(G)$.

It is known that for every continuous function $G: X \to \mathbb R$ there exists at least one equilibrium state.  Under the additional assumption that $G$ is H\"older continuous then there is a unique measure $\mu_G \in \mathcal M_\phi(X)$ satisfying (1.4) which is the Gibbs measure for $G$.

The SRB measure for a hyperbolic attracting  flow  $\phi_t: X \to X$
is  a Gibbs measure associated to the  H\"older continuous function $\Phi: X \to \mathbb R$ defined by  
$$\Phi(x) = - \lim_{t \to 0^+} \frac{1}{t}\log \left|\det (D\phi_t|E_x^u)\right|,$$
which is  again  called the unstable expansion coefficient, where $E^u_x$ is the unstable bundle for the flow.
The analogue of Theorem \ref{thm:sinai} for flows was proved by Bowen and Ruelle \cite{BR}.
 
  More generally, given any  H\"older continuous function $G : X \to \mathbb R$  we give
a construction of the  Gibbs measure $\mu_G$ which is analogous to that for diffeomorphisms in 
 Theorem  \ref{thm:gibbsdiffeo}:  


\begin{thm}\label{thm:gibbsflow}
Let $\phi_t:  X \to X$ be a $C^{1+\alpha}$ mixing hyperbolic attracting flow  and  let  $G: X \to \mathbb R$
be a H\"older continuous function. 
Given $x\in X$ and $\delta > 0$ 
consider 
the  family of probability measures $(\lambda_T)_{T>0}$ 
supported on a local unstable manifold $W^u_\delta(x)$ 
and absolutely continuous with respect to the induced volume $\lambda = \lambda_{W^u_\delta(x)}$
with densities
$$
\frac{d\lambda_T}{d\lambda}(y) :=  
\frac{
 \exp\left( \int_0^T (G - \Phi)(\phi_vy))dv \right) 
 }{
 \int_{W_\delta^u(x)}
 \exp\left( \int_0^T (G - \Phi)(\phi_vy))dv \right)
 d\lambda_{W_\delta^u(x)}(z)
 }
 \hbox{ for } y \in  W^u_\delta(x).
$$
Then the averages  
$$\mu_T :=\frac{1}{T}\int_0^T (\phi_t)^*\lambda_T dt, \quad T > 0,
$$
 converge in the weak star topology to $\mu_G$ as $T \to +\infty$.
\end{thm}

As in the case of hyperbolic  diffeomorphisms each of the 
absolutely continuous reference measures $\lambda_T$ depends on $T>0$.
More generally  if we only assume that $G: X \to \mathbb R$ is continuous then the limit points of $\mu_T$ are equilibrium states for $G$.
Similarly, the contribution  of the  term involving $\Phi$ comes from the change of variables for the time $T$ flow diffeomorphism 
 $\phi_T: W^u_\delta(x) \to \phi_T(W^u_\delta(x))$ for which the Jacobian is 
$$\left|\det (D\phi_T|E_y^u)\right|
=  
 \exp\left( -\int_0^T \Phi(\phi_vy) dv\right)
 $$
  and  we can rewrite
 $$
 \lambda_T(A) = 
 \frac{\int_{\phi_T(A)}  \exp\left( \int_0^T G(\phi_{-v}y) dv\right)d\lambda_{\phi_TW_\delta^u(x)}(y)}
 {\int_{\phi_T(W_\delta^u(x))}{ \exp\left( \int_0^T G(\phi_{-v}z) dv\right)d\lambda_{\phi_TW_\delta^u(x)}(z)}}
 \hbox{  for Borel  $A \subset W_\delta^u(x)$.}
 $$ 

\begin{remark}
In Theorem   \ref{thm:gibbsflow}
it is possible to replace the local unstable manifold $W^u_\delta(x)$ by an embedded
 curve $\mathcal C \subset X$ of dimension $\dim E^u$ provided it does not 
 wholly 
 lie inside a single stable manifold. 
 This can again be done by uniformly approximating $\mathcal C$
 by pieces of unstable manifold.
 \end{remark}
 

Whereas the original proof(s) of Theorem \ref{thm:sinai} used symbolic dynamics, we use a more geometric approach.  The main ingredients are estimates on the rate of growth of the volume of the 
local unstable manifolds (i.e., $\hbox{\rm vol}(f^nW_\delta^u(x))$ as $n \to +\infty$ for hyperbolic attracting diffeomorphisms and $\hbox{\rm vol}(\phi_T W_\delta^u(x))$ as $T \to +\infty$ for hyperbolic attracting flows) and entropy bounds adapted from work of Misiurewicz.
In particular, this gives a new proof of Theorem \ref{thm:sinai}.  

\section{SRB measures and measures  of maximal entropy}
In this section we will concentrate on  two  particularly important examples of Gibbs measures   
for each of  the cases of diffeomorphisms and flows.   These are the SRB measures
(which we already encountered in Theorem \ref{thm:sinai}) and the measure of maximal entropy, 
also known as the Bowen-Margulis measure.

\subsection{SRB measures and measures  of maximal entropy  for hyperbolic diffeomorphisms}

The importance of the SRB measure is that it  describes the behaviour of orbits for typical points in the basin of 
attraction \cite{ruelle}. 
On the other hand, the  measure of maximal entropy  describes the distribution of various quantities which grow at a rate corresponding to the topological entropy, such as the periodic points  \cite{bowen-periodic},  \cite{parry1}.

\begin{example}[SRB-measures]
In the special case that the Gibbs measure is  the SRB measure we can take 
 $G = \Phi = -  \log \left|\det ( Df|E_x^u)\right|$ and thus the weights reduce to 
$$
\frac{d\lambda_n}{d\lambda}(y) =  \exp\left( \sum_{i=0}^{n-1} (G - \Phi)(f^iy)\right)
 =1
\hbox{ for all $y \in W^u_\delta(x)$}
 $$
and then by definition
$\lambda_n = \lambda$, for all $n \geq 1$.  In particular, 
$$
\mu_n = 
\frac{1}{n}\sum_{k=0}^{n-1} f_*^k \lambda_n = \frac{1}{n}\sum_{k=0}^{n-1} f_*^k \lambda, 
$$
for $n \geq 1$,  
converges to $\mu_{SRB}$, 
and so
Theorem \ref{thm:gibbsdiffeo} reduces to Theorem \ref{thm:sinai}.
\end{example}


\begin{example}[Measure of Maximal Entropy]
In  the special case that the Gibbs measure is the measure of maximal entropy  
we can take  $G = 0$.  In particular, (1.1) now reduces to  $P(0) = h(f)$, the topological entropy.
Furthermore,  the sequence of densities  $(\lambda_n)_{n=1}^\infty$
 is  given by
$$
\begin{aligned}
\frac{d\lambda_n}{d\lambda}(y) &:=  
\frac{
 \exp\left( -\sum_{i=0}^{n-1} \textcolor{black}{\Phi}(f^iy) \right) 
 }{
 \int_{W_\delta^u(x)}
 \exp\left( -\sum_{i=0}^{n-1} \textcolor{black}{\Phi}(f^iz) \right) 
 d\lambda_{W_\delta^u(x)}(z)
 }\cr
&=
\frac{
\left|\det ( Df^n|E_y^u)\right|
 }{
 \int_{W_\delta^u(x)}
\left|\det (Df^n|E_z^u)\right|
 d\lambda_{W_\delta^u(x)}(z)
 }
 \hbox{ for  $y \in  W^u_\delta(x)$.}
\cr
\end{aligned}
$$
In particular,  the averages 
$$
\mu_n = 
\frac{1}{n}\sum_{k=0}^{n-1} f_*^k
\left(\frac{
\left|\det ( Df^n|E_x^u)\right|
 }{
 \int_{W_\delta^u(x)}
\left|\det (Df^n|E_x^u)\right|
  d\lambda_{W_\delta^u(x)}(x)
 }  \lambda \right), \quad n \geq 1, 
$$ converge  to the measure of maximal entropy $\mu_{BM}$.
\end{example}

For typical hyperbolic attracting diffeomorphisms the measure of maximal entropy  and the  SRB measure are singular with respect to each other.

\begin{example}[Toral automorphism]
In the particular case that $M = \mathbb T^2$ then any Anosov diffeomorphism   
$f:\mathbb T^2  \to \mathbb T^2$
is topologically conjugate to a linear hyperbolic diffeomorphism.  In the special 
case  that  $f$ is $C^1$ conjugate 
to a   linear hyperbolic toral automorphism then the  function $\Phi: \mathbb T^2 \to \mathbb R$ is 
{\color{black} (cohomologous to a)}
constant and the measure of maximal entropy coincides with the SRB measure and is equal to the normalized Haar measure.  Otherwise,  the measure of maximal  entropy  and the  SRB measure are singular with respect to each other.
\end{example}

\subsection{SRB measures and measures  of maximal entropy  for hyperbolic attracting  flows}

As in the case of hyperbolic attracting  diffeomorphisms we can  consider  the two special cases of Gibbs measures for hyperbolic attracting flows, i.e.,  the SRB measure and the measure of maximal entropy.

\begin{example}[SRB-measures]
In the special case that the Gibbs measure is  the SRB measure we can take 
 $G(x) = \Phi(x)  = -  \lim_{t\to 0} \frac{1}{t}\log \left|\det ( D\phi_t|E_x^u)\right|$ and thus the densities $(\lambda_T)_{T>0}$  reduce to 
$$
\frac{d\lambda_T}{d\lambda}(y) =  \exp\left(\int_0^T (G - \Phi)(\phi_uy) du\right)
 =1
\hbox{ for all $y \in W^u_\delta(x)$}
 $$
and then by definition
$\lambda_T = \lambda$, for all $T>0$.  In particular, 
$$
\mu_T = 
\frac{1}{T} \int_0^T (\phi_{\textcolor{black}{t}})^{\textcolor{black}{*}} \lambda_{\textcolor{black}{T}} dt
=\frac{1}{T} \int_0^T (\phi_t)^{\textcolor{black}{*}} \lambda dt
$$
for $T > 0$ 
converges to 
$\mu_{SRB}$
and 
Theorem \ref{thm:gibbsflow} reduces to the analogue of  Theorem 
 \ref{thm:sinai}
for flows, which was originally proved by \cite{BR}. 
\end{example}

The second example is the measure of maximal entropy, or Bowen-Margulis measure.

\begin{example}[Measure of Maximal Entropy]
In  the special case that the Gibbs measure is the measure of maximal entropy  
we can take  $G = 0$.  Thus  the reference  measures $(\lambda_T)_{T>0}$
given by
$$
\begin{aligned}
\frac{d\lambda_T}{d\lambda}(y) &:=  
\frac{
 \exp\left( -\int_0^T\Phi(\phi_u y)du \right) 
 }{
 \int_{W_\delta^u(x)}
 \exp\left( -\int_0^T\Phi(\phi_u y)du \right) 
 d\lambda_{W_\delta^u(x)}(z)
 }\cr
&=
\frac{
\left|\det ( D\phi_T|E_y^u)\right|
 }{
 \int_{W_\delta^u(x)}
\left|\det (D\phi_T|E_z^u)\right|
 d\lambda_{W_\delta^u(x)}(z)
 }
 \hbox{ for  $y \in  W^u_\delta(x)$.}
\cr
\end{aligned}
$$
In particular, 
$$
\mu_T = 
\frac{1}{T}\int_0^T  (\phi_t)^*
\left(\frac{
\left|\det ( D\phi_{\textcolor{black}{T}}|E_y^u)\right|
 }{
 \int_{W_\delta^u(x)}
\left|\det (D\phi_{\textcolor{black}{T}}|E_z^u)\right|
 d\lambda_{W_\delta^u(x)}(z)
 }  \lambda \right) dt, \quad T > 0,
$$ converges to the measure of maximal entropy $\mu_{BM}$.
\end{example}

\begin{example}[Geodesic flows]
  The classical example of an attracting hyperbolic  flow 
  (or even an Anosov flow)
  is the geodesic flow on a 
  compact surface  with negative  curvatures \cite{anosov}.
  In the special case of the  curvatures  being constant the function 
   $\Phi$ is constant and the measure of maximal entropy coincides with the SRB measure and is equal to the normalized volume.  
 For surfaces of negative non-constant curvature  the SRB measure is still equivalent to the volume but the measure of maximal  entropy is singular with respect to the volume.
 \end{example}

\section{Definitions  and simple  examples}
We now recall  the formal definitions of  hyperbolic attracting  diffeomorphisms and Anosov diffeomorphisms. 
  We then give some simple examples which illustrate Theorems
\ref{thm:gibbsdiffeo} and \ref{thm:gibbsflow}

\subsection{Definition of hyperbolic diffeomorphisms}


We   recall the general definition of a  hyperbolic attractor.
Let $f: M \to M$ be a $C^{1+\alpha}$ diffeomorphism on a compact \textcolor{black}{Riemannian} manifold, and let $X \subset M$ be a closed $f$-invariant set.
\begin{definition}
The map 
 $f:X \to X$ is called a  {\it \textcolor{black}{mixing} hyperbolic attracting diffeomorphism} if: 
\begin{enumerate}
\item 
there exists a continuous splitting $T_X M = E^s\oplus E^u$ and $C > 0$ and $0< \lambda < 1$ such that 
$$\|Df^n|E^s\| \leq C \lambda^n \hbox{ and }\|Df^{-n}|E^u\| \leq C \lambda^n$$ for $n \geq 0$;
\item there exists an open set $X \subset U \subset M$ such that 
$X = \cap_{n=0}^\infty f^{n}U$;
\item
$f: X \to X$ is topologically mixing; and
\item
the periodic orbits for 
$f: X \to  X$ are dense in $X$.
\end{enumerate}
\end{definition}

In the particular case $X=M$  the diffeomorphism  $f$ is an  Anosov diffeomorphism:
\begin{definition}
A mixing diffeomorphism
 $f:M \to M$ is called  {\it Anosov} if: 
\begin{enumerate}
\item 
there exists a continuous splitting $TM = E^s\oplus E^u$ and $C > 0$ and $0< \lambda < 1$ such that 
$$\|Df^n|E^s\| \leq C \lambda^n \hbox{ and }\|Df^{-n}|E^u\| \leq C \lambda^n$$ for $n \geq 0$; and 
\item
$f: X \to X$ is topologically mixing.
\end{enumerate}
\end{definition}


We next recall some  classic examples of hyperbolic diffeomorphisms.

\begin{example}[Axiom A]
An Axiom A diffeomorphism is a diffeomorphism 
 for which the non-wandering set  $\Omega$ is a finite union of hyperbolic sets 
and attracting or repelling periodic points \cite{smale}. 
\end{example}  

A more concrete example is the following.

\begin{example}[Solenoid]
One of the simplest examples of a hyperbolic  attractor (which is not Anosov) is the solenoid
constructed in terms of a solid torus mapped inside itself (whose interior plays the role of the neighbourhood $U$), {\color{black} see \cite{williams} and   \cite{kh}, \S 7.1. More concretely,  let $M = \{(\theta, x,y) \in \mathbb R/\mathbb Z \hbox{ : } x^2 + y^2 \leq 1 \}$be a solid torus and let $f: M \to M$ be defined by
$$
f(\theta, x, y) =
\left(
2\theta, \frac{x}{10} + \frac{\cos \theta}{2}, \frac{x}{10} + \frac{\sin \theta}{2}
\right).
$$
It is easily seen that $f(M) \subset \hbox{\rm int}(M)$ and the $f$-invariant set
 $\Lambda = \cap_{n=0}^\infty f^n\hbox{\rm int}(M)$ is an attractor (see \cite{kh}, Proposition 17.1.2).
In this example the unstable manifold is one dimensional and the stable manifold is two dimensional.
The SRB induces the usual Haar measure on the unstable manifolds, which are locally parameterized by the $\theta$-coordinate.  In this case the measure of maximal entropy coincides with the SRB-measure.  
}

\end{example}

\subsection{\textcolor{black}{Example} illustrating Theorem \ref{thm:gibbsdiffeo}}
It is illuminating to consider \textcolor{black}{the following} simple \textcolor{black}{example}.
 
 \begin{example}[Arnol'd  CAT map]
 Let $f: \mathbb T^2 \to \mathbb T^2$ be a linear hyperbolic toral automorphism of the form
 $
 f\left( (x,y) + \mathbb Z^2 \right) =  (ax + by, cx + dy) + \mathbb Z^2 
 $
 where $a,b,c,d \in \mathbb Z$ with $ad-bc = 1$ and  $|a+d| \neq 2$.
 The unstable manifolds are line segments of irrational slope $\alpha$, say,  and $\Phi$ is a constant function.
 Thus, for any H\"older function $G: \mathbb T^2 \to \mathbb R$ we can write
 $$
\frac{d\lambda_n}{d\lambda}(y) :=  
\frac{
 \exp\left( \sum_{i=0}^{n-1} G(f^iy)) \right) 
 }{
 \int_{W_\delta^u(x)}
 \exp\left( \sum_{i=0}^{n-1} G(f^iz) \right) 
 d\lambda(z)
 }
 \hbox{ for } y \in  W^u_\delta(x),\eqno(3.1)
$$
for $n \geq 1$.
If we take the specific choice  that 
$$
W_\delta^u(x)
= \left\{
t (\ell, \alpha\ell) \hbox{ : } 0 \leq t \leq 1
\right\}  + \mathbb Z^2
$$
 is a line segment  from $(0,0)$ to $(\ell, \alpha\ell) + \mathbb Z^2$, for some $\ell > 0$,  then we see that
  $$
  f^nW_\delta^u(x)
  =  \left\{
t(e^{hn} \ell, e^{hn} \alpha\ell)   \hbox{ : } 0 \leq t \leq 1
\right\}  + \mathbb Z^2
  $$ is the line segment from  $(0,0)$ to $(e^{hn} \ell, e^{hn} \alpha\ell) + \mathbb Z^2$, where $h>0$ is the topological entropy of $f$.  

 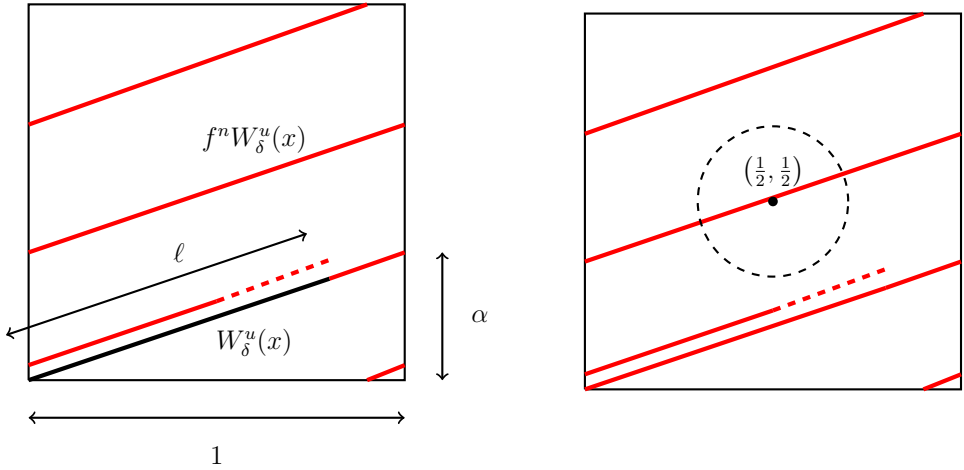
\begin{figure}[h!]
\centerline{
\begin{tikzpicture}[thick,scale=1.0, every node/.style={scale=0.85}]
\draw (0,0) -- (5,0)--(5,5)--(0,5)--(0,0);
\draw[ultra thick] (0,0) -- (4,1.35);
\draw[ultra thick,red] (4,1.35) -- (5,1.7);
\draw[ultra thick, red] (0,1.7) -- (5,3.4);
\draw[ultra thick, red] (0,3.4)--(4.5, 5);
\draw[ultra thick, red] (4.5, 0)--(5,0.2);
\draw[ultra thick, red] (0,0.2)-- (2.5,1.05);
\draw[ultra thick, red, dashed] (2.5,1.05)-- (4,1.6);
\draw[<->] (0,-0.5) -- (5,-0.5) ;
\draw[<->] (5.5,0) -- (5.5,1.7) ;
\node at (2.5,-1) {$1$};
\node at (6,0.85) {$\alpha$};
\node at (3,0.5) {$W^u_\delta(x)$};
\node at (3,3.2) {$f^nW^u_\delta(x)$};
\draw[<->] (-0.3,0.6) -- (3.7,1.95);
\node at (2,1.7) {$\ell$};
\end{tikzpicture}
\hskip 1cm
\begin{tikzpicture}[thick,scale=1.0, every node/.style={scale=0.8}]
\draw (0,0) -- (5,0)--(5,5)--(0,5)--(0,0);
\draw[ultra thick,red] (0,0) -- (4,1.35);
\draw[ultra thick,red] (4,1.35) -- (5,1.7);
\draw[ultra thick, red] (0,1.7) -- (5,3.4);
\draw[ultra thick, red] (0,3.4)--(4.5, 5);
\draw[ultra thick, red] (4.5, 0)--(5,0.2);
\draw[ultra thick, red] (0,0.2)-- (2.5,1.05);
\draw[ultra thick, red, dashed] (2.5,1.05)-- (4,1.6);
\draw[dashed] (2.5,2.5) circle (1cm);
\draw[fill] (2.5,2.5) circle (0.05cm);
\node at (2.5,2.9) {$\left(\frac{1}{2}, \frac{1}{2} \right)$};
\node at (2.5,-1) {};
\end{tikzpicture}
}
\caption{A piece of local unstable manifold $W^u_\delta(x)$ of length $\ell$
for the Arnol'd CAT map
 and its image 
$f^nW^u_\delta(x)$}
\end{figure}

In the case that $G=0$ the SRB measure and the measure of maximal entropy coincide, and are both equal to the Haar measure.
We can consider two balls:
\begin{enumerate}
\item
$B_1 = B\left((0,0), \frac{1}{3} \right)$ of radius $\epsilon = \frac{1}{3}$ centred at $(0,0)$; and 
\item 
$B_2 = B\left((\frac{1}{2}, \frac{1}{2}), \frac{1}{3} \right)$
 of radius $\epsilon = \frac{1}{3}$
 centred at $(\frac{1}{2}, \frac{1}{2})$, 
\end{enumerate}
which both have $\lambda$-measure $\pi^2/9$.  Figure 3 shows the values of $\mu_n(B_1)$ and $\mu_n(B_2)$ for each of the two potentials 
$G(x,y)=0$ and $G(x,y) = \frac{1}{10} \sin \left( 2\pi x \right)$.


\begin{figure}[h!] 
\centerline{
\includegraphics[width=6cm]{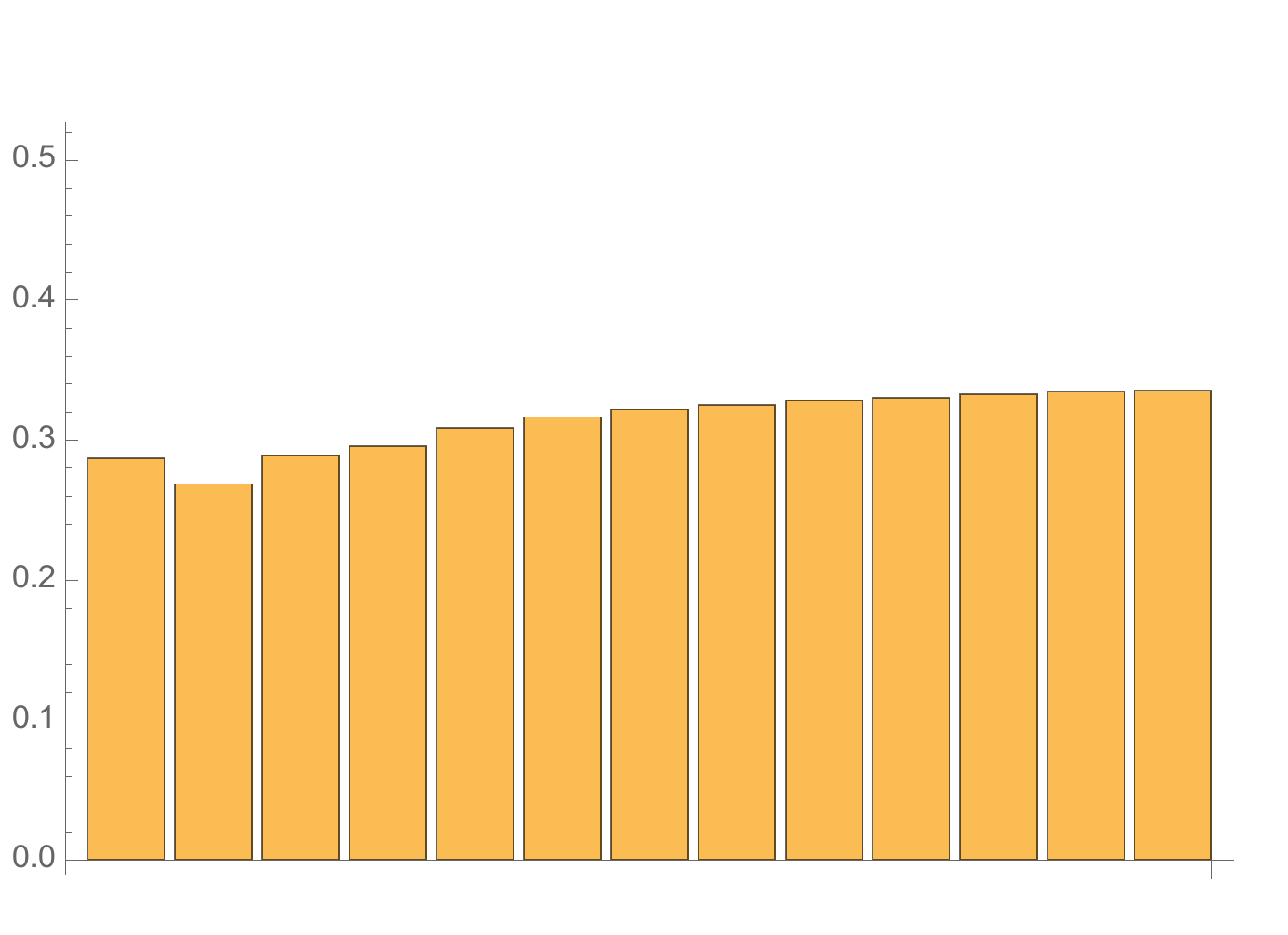}
\hskip 1cm
\includegraphics[width=6cm]{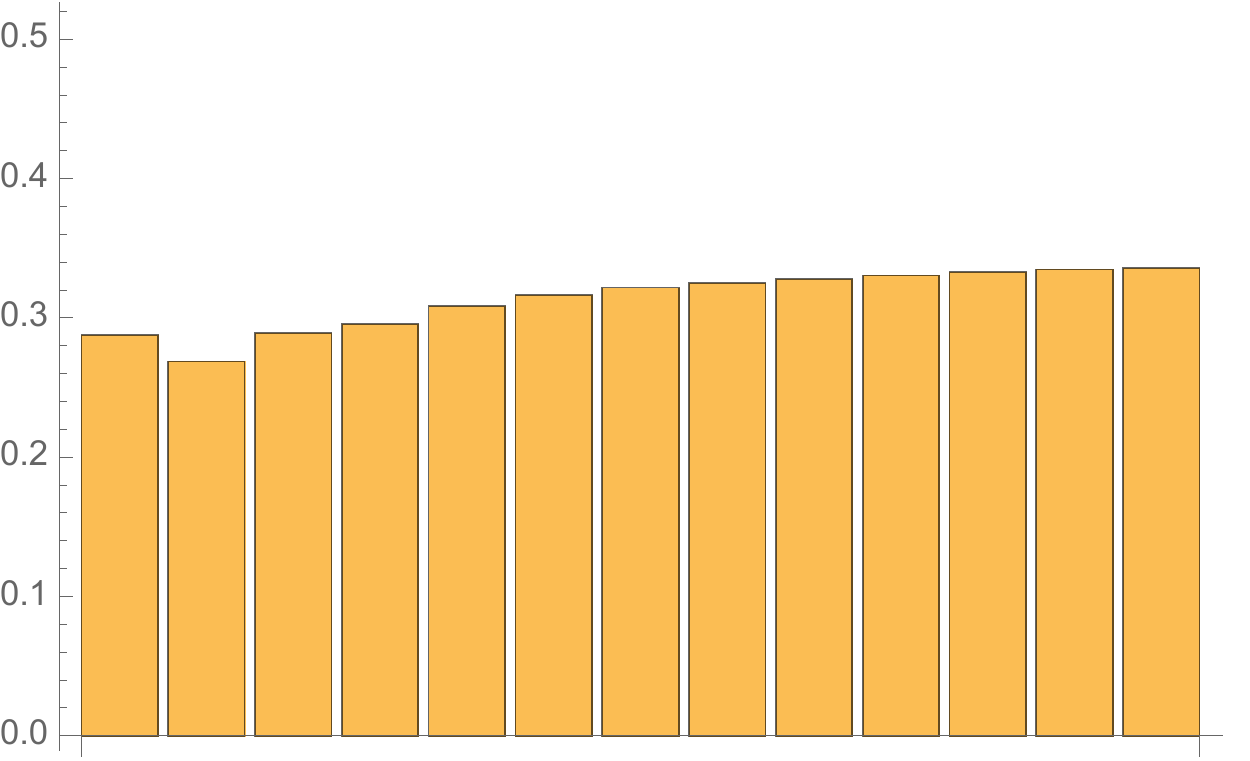}
}
\vskip 0.75cm
\centerline{
\includegraphics[width=6cm]{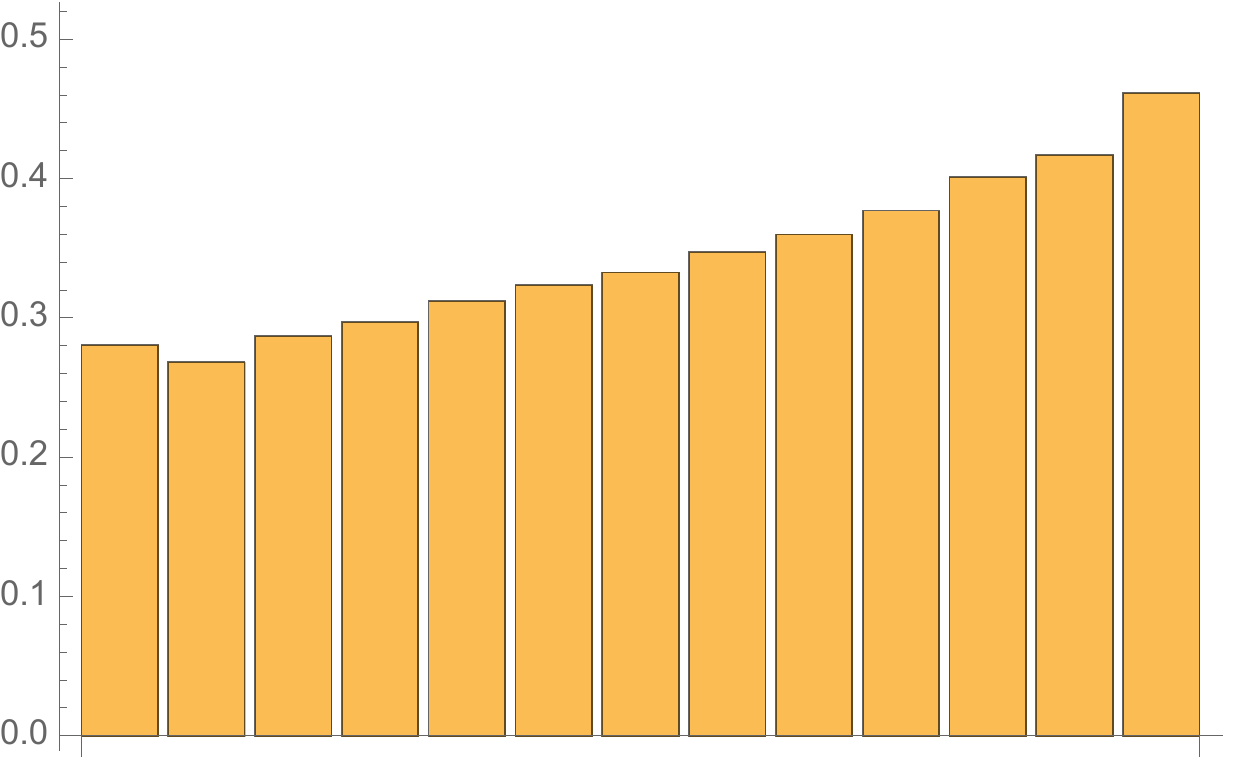}
\hskip 1cm 
\includegraphics[width=6cm]{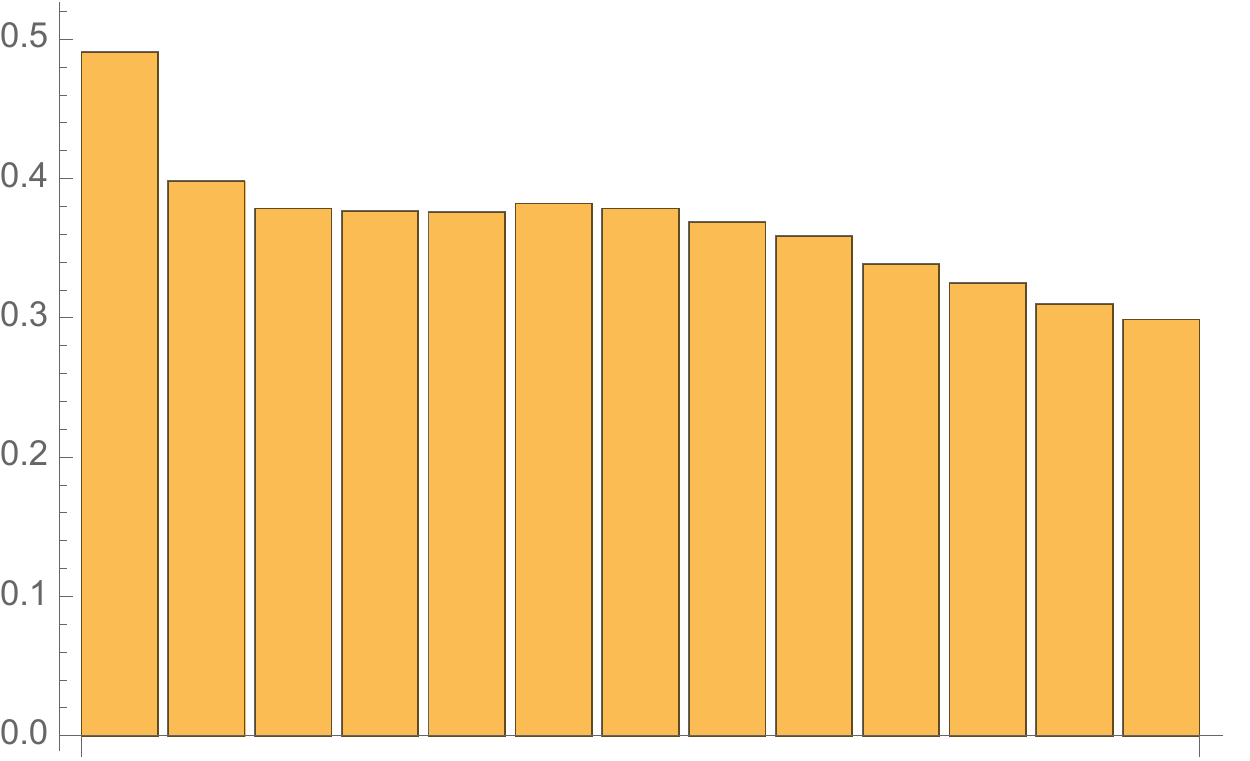}
}
\caption{(a) When $G=0$ the two top bar charts represent  $\mu_n(B_1)$ and   $\mu_n(B_2)$  for $n=1, \cdots, 12$ and the convergence to $\frac{\pi^2}{9}$ can be seen.
(b) When $G(x,y) =  \frac{1}{10}\sin \left( 2\pi x \right)$ the two  lower bar charts represent  $\mu_n(B_1)$ and   $\mu_n(B_2)$ for $n=1, \cdots, 12$.
{\color{black}
The sequences $\mu_n(B_1)$ and $\mu_n(B_2)$ converge to the measures of the balls $B_1$ and $B_2$ with respect to the corresponding Gibbs measure.  It its empirically clear these values are different,
 but the precise values are unknown to us.}
}
\end{figure}




 \end{example}

\subsection{Example illustrating Theorem  \ref{thm:gibbsflow}}

We return to the important example of geodesic flows.

  \begin{example}[Geodesic flow]
Let $\phi_t: M \to M$ be the geodesic flow on the
three dimensional 
 unit tangent bundle 
$M=SV$
of 
a compact surface  $V$ of 
negative curvature.  
In this case the unstable manifolds correspond to horocycles.
Given $x \in V$ we let $\mathcal C = S_xV$ be the fibre above $x \in V$.
We define a family of   densities $(\lambda_T)_{T>0}$ on $\mathcal C$ by 
$$
 \lambda_T(A) = 
 \frac{\int_{\phi_T(A)}  \exp\left( \int_0^T \textcolor{black}{G}(\phi_{\textcolor{black}{-v}}x) dv\right)d\lambda_{\mathcal C}(x)}
 {\int_{\phi_T(\mathcal C)}{ \exp\left( \int_0^T \textcolor{black}{G}(\phi_{\textcolor{black}{-v}}x) dv\right)d\lambda_{\mathcal C}}(x)}
 \hbox{  for Borel  $A \subset \mathcal C$.}
 $$ 
Then
by  Theorem  \ref{thm:gibbsflow}
 the averages of the push forwards 
$$\mu_T :=\frac{1}{T}\int_0^T (\phi_t)^{\textcolor{black}{*}}\lambda_T dt,    \quad  T > 0, 
$$
 converge in the weak star topology to $\mu_G$ as $T \to +\infty$.

\begin{enumerate}
\item If $G=0$ then the limiting measure $\mu_G$ is the Bowen-Margulis measure.
\item If $G=\Phi$ then the limiting measure $\mu_G$ is the Liouville measure (or, equivalently, the SRB measure).
\end{enumerate}
If $V$ is a surface of constant curvature $\kappa=-1$ then 
$\Phi = 1$ is a constant function and the Bowen-Margulis measure is equal to the Liouville measure. 

\end{example}

\section{Preliminaries for Theorem \ref{thm:gibbsdiffeo}}
In this section we  collect together some  results  for hyperbolic attracting diffeomorpisms in anticipation of \textcolor{black}{the proof of} Theorem \ref{thm:gibbsdiffeo} in the next section.

\subsection{Stable and unstable manifolds}  
We next recall the definition of the stable and unstable manifolds. 
  
 \begin{definition}
We associate to each $x \in X$ the  {\it stable manifold} defined by
$$
W^{s}(x) = \{y \in M \hbox{ : } d(f^nx, f^ny) \to 0 \hbox{ as } n \to +\infty\}
$$
and  {\it unstable manifold}  defined by
$$
W^{u}(x) = \{y \in M \hbox{ : } d(f^{-n}x, f^{-n}y) \to 0  \hbox{ as } n \to +\infty\}.
$$
For sufficiently small $\delta > 0$ we define local versions 
$W_\delta^{s}(x) \subset W^{s}(x)$ and  $W_\delta^{u}(x) \subset W^{u}(x)$
defined by 
$$
W_\delta^{s}(x) = \{y \in M \hbox{ : } d(f^nx, f^ny) \leq \delta, \forall n \geq 0\}
$$
and 
$$
W_\delta^{u}(x) = \{y \in M \hbox{ : } d(f^{-n}x, f^{-n}y) \leq \delta, \forall n \geq 0\}.
$$
\end{definition}

The basic properties of the (local) stable and unstable manifolds are the following.

\begin{lemma}[Hirsch-Pugh \cite{hp}]\label{stab}
For each $x \in X$ and sufficiently small $\delta > 0$, the sets:
\begin{enumerate}
\item
$W^{s}(x)$ and $W^{u}(x)$
are $C^{1}$ immersed submanifolds of dimension 
$\dim (E^s)$ and $\dim (E^u)$, respectively, with $T_xW^s = E_x^s$
and $T_xW^u = E_x^u$; and 
\item
$W_\delta^s(x)$ and $W_\delta^{u}(x)$
are $C^{1}$ embedded  disks  of dimension 
$\dim (E^s)$ and $\dim (E^u)$, respectively, with $T_xW_\delta^s(x) = E_x^s$
and $T_xW^u_\delta(x) = E_x^u$.
\end{enumerate}
\end{lemma}


\subsection{Entropy}
We now describe  some results on the entropy of invariant probability measures.  For the duration of this subsection we will only require 
the weaker assumption that $f: X \to X$ is a homeomorphism.

We begin with some standard   definitions.

\begin{definition}
Given a finite measurable partition $\mathcal P = \{P_1, \cdots, P_k\}$ and a probablity measure $\nu$ we can associate the
 entropy of the partition defined by 
 $$
 H_\nu(\mathcal P)
 =  - \sum_{i=1}^k  \nu(P_i)  \log \nu(P_i). 
 $$
 \end{definition}
 
 \noindent
 We let $\bigvee_{i=0}^{n-1} f^{-i} \mathcal P = \{P_{i_0} \cap f^{-1} P_{i_1} \cap \cdots  \cap f^{-(n-1)} P_{i_{n-1}}
 \hbox{ : } 1 \leq  i_0, \cdots, i_{k-1} \leq k
 \}$ be the refinement of the partitions $\mathcal P, f^{-1}\mathcal P, \cdots, f^{-(n-1)}\mathcal P$.

 \begin{definition}
We can then define the entropy associated to the partition $\mathcal P$ by 
 $$
 h_\nu(\mathcal P) = \lim_{n \to +\infty} \frac{1}{n} H_\nu\left( \bigvee_{i=0}^{n-1} f^{-i} \mathcal P\right)
 $$
Lastly, the entropy with respect to the measure is defined by
$$
h(\nu) = \sup_{\mathcal{P}}\{h_{\nu}(\mathcal{P}) \hbox{ : $\mathcal{P}$ is a countable partition with $H_{\nu}(\mathcal{P}) < \infty$}\}.
$$
In the case that $\mathcal P$ is a generating partition we have that $h(\nu) = h_\nu(\mathcal P)$ is the entropy of the measure $\nu$.
 \end{definition}
 
 {\color{black}
 
 The following lemma will be useful later.
 
  \begin{lemma}
  $$
  q H_{\lambda_n} \left(
  \bigvee_{h=0}^{n-1} f^{-h} \mathcal P 
 \right)
 \leq 
n   H_{\mu_n}\left(
  \bigvee_{i=0}^{q-1} f^{-i} \mathcal P 
 \right)
 +
 2q^2 {\color{black} \log } \hbox{\rm Card} (\mathcal P).
  $$
  \end{lemma}
}

\begin{proof}
Following  a construction  of Misiurewicz (see \cite{walters}, p.220) we proceed as follows.
Given $0 < q < n$ choose  $0 \leq k < q $  and then let $a = a(k)  = [(n-k)/q]$.  We can then partition 
$\{0,1,2, \cdots, n-1\}$  by 
 $$
 \{0, \cdots, k\} \cup
 \bigcup_{l=0}^{a-1} \{ lq + k, \ldots,  (l+1)q + k-1 \}
  \cup\{ \textcolor{black}{aq + k}, \ldots, n-1\}.
 $$
 We can use this to rewrite the refinement as 
 $$
 \bigvee_{h=0}^{n-1} f^{-h} \mathcal P 
 =
  \bigvee_{t=0}^{k} f^{-t} \mathcal P \vee
  \left( 
  \bigvee_{r=0}^{a-1}
  f^{-(rq+k)}
  \left(\bigvee_{i=0}^{q-1} f^{-i}\mathcal P\right)
   \right)
 \vee
 \bigvee_{u= \textcolor{black}{a}q + k}^{  n-1} f^{-u} \mathcal P.
 $$
Thus we can write
 $$
 \begin{aligned}
 H_\nu\left(
  \bigvee_{h=0}^{n-1} f^{-h} \mathcal P 
 \right)
& \leq 
  H_\nu\left( \bigvee_{t=0}^{k} f^{-t} \mathcal P
  \right)
  +
 \sum_{r=0}^{a-1} H_\nu\left( 
  f^{-(rq+k)}
  \left(\bigvee_{i=0}^{q-1} f^{-i} \mathcal P\right)
   \right)\cr
&\qquad  +
   H_\nu\left(
    \bigvee_{u= \textcolor{black}{a}q + k}^{ n-1} f^{-u} \mathcal P
   \right)
   \end{aligned}
   \eqno(4.1)
 $$
 (cf. Walters \cite{walters})
 as $k < q$ and $aq + k \geq n - q$ we have the bounds
 $$
   H_\nu\left( \bigvee_{t=0}^{k} f^{-t} \mathcal P
  \right), 
   H_\nu\left(
    \bigvee_{u= \textcolor{black}{a}q + k}^{  n-1} f^{-u} \mathcal P
   \right) \leq q {\color{black} \log } \hbox{\rm Card} (\mathcal P).
 $$
Summing both sides of (4.1) over $k$, we  can write
 $$
 \begin{aligned}
 q H_\nu\left(
  \bigvee_{h=0}^{n-1} f^{-h} \mathcal P 
 \right)
& \leq 
 \sum_{k=0}^{q-1}
 \left(
 \sum_{r=0}^{a-1} H_\nu\left( 
  f^{-(rq+k)}
  \left(\bigvee_{i=0}^{q-1} f^{-i}\mathcal P\right)
   \right)
 + 
 2 q  {\color{black} \log } \hbox{\rm Card} (\mathcal P)  \right)\cr
 &\leq \sum_{m=0}^{n-1}
H_{\textcolor{black}{f^m_*}\nu}
  \left(\bigvee_{i=0}^{q-1} f^{-i}\mathcal P\right)
 + 
 2q^2  {\color{black} \log } \hbox{\rm Card} (\mathcal P),
 \end{aligned}\eqno(4.2)
  $$
  {\color{black} using \cite{walters}, Remark 2 (iii) \S 8.2   and that 
  $H_\nu(f^{-s}\xi) = H_{f_*^{s}\nu}(\xi)$ ($s \geq 0$) for any partition $\xi$}.
  If we let $\nu = \lambda_n$ and $\mu_n: = \frac{1}{n} \sum_{\textcolor{black}{m=0}}^{n-1} f_*^{\textcolor{black}{m}}\lambda_n$ then by concavity of the entropy
  $$
  H_{\mu_n}
  \left(\bigvee_{i=0}^{q-1} f^{-i}\mathcal P\right)
  \geq  \frac{1}{n} \sum_{\textcolor{black}{m}=0}^{n-1} 
    H_{f_*^{\textcolor{black}{m}}\lambda_n}
  \left(\bigvee_{i=0}^{q-1} f^{-i}\mathcal P\right).
  \eqno(4.3)
  $$
    (cf. \cite{walters}, Remark 1, \S8.2).
{\color{black}  Comparing (4.2) and (4.3) concludes the proof.}
   \end{proof}
   

\subsection{Pressure and growth}
The pressure $P(G)$ has various alternative interpretations in terms of the growth of different  quantities.
We will need the following  particular variant on these.

\begin{proposition}\label{pressure}
Let $f: X \to X$ be a mixing hyperbolic attracting diffeomorphism. 
For any  continuous function $G: X \to \mathbb R$
we have 
$$
P(G)
= \lim_{n \to +\infty} \frac{1}{n}\log  \int_{W_\delta^{u}(x)}
 \exp\left( \sum_{k=0}^{n-1} (G-\Phi)(f^ky) \right) 
 d\lambda_{W_\delta^u(x)}(y). \eqno(4.4)$$
\end{proposition}

Again, we can  use the change of variables formula  to rewrite this in the equivalent form
$$
P(\textcolor{black}{G})
= \lim_{n \to +\infty} \frac{1}{n}\log  \int_{f^nW_\delta^{u}(x)}
\exp\left( \sum_{k=1}^{n} G(f^{-k}y) \right)
  d\lambda_{f^nW_\delta^u(x)}(y).
 $$
 
 \begin{proof}[Proof of Proposition \ref{pressure}]
 We begin with the following standard result.

\begin{lemma}\label{cover0}
\textcolor{black}{F}or any $\epsilon > 0$ there exists \textcolor{black}{an} $m > 0$ such that $f^mW^u_\delta(x)$ 
is $\epsilon$-dense in $X$.  In particular, we can assume that
 $X = \cup_{y \in f^mW^u_\delta(x)} W_\epsilon^s(y).$
\end{lemma}

\begin{proof}
This is a consequence of the minimality of the unstable foliation for mixing hyperbolic attracting diffeomorphisms and the local product structure.  
\end{proof}

We recall that the pressure $P(G)$ can also  be written in terms of the growth rates of 
spanning sets and separated sets. Recall that given $\epsilon > 0$ and 
$n \geq 1$  an
 $(n,\epsilon)$-spanning set $S \subset X$ 
is such that $ \cup_{x\in S} B(x,n,\epsilon)$ covers $X$,  where 
$B(x,n,\epsilon) := \cap_{k=0}^{n-1}f^{-k} B(f^kx, \epsilon)$  is a   Bowen ball. 
On the other hand  a\textcolor{black}{n} $(n,\epsilon)$-separated  set $\Sigma \subset X$ 
is such that \textcolor{black}{$d_n(x,y) > \epsilon$} for $x,y \in \Sigma$
 (and, in particular, $B(x,n,\epsilon/2)$, $x \in \Sigma$, are disjoint
in $X$).

\begin{lemma} 
Given 
$n \geq 1$ and $\epsilon > 0$,  let 
$$
Z_0(n, \epsilon) = \inf\left\{
\sum_{y \in S} \exp(G^n(y))
\hbox{ : } S \hbox{ is an $(n, \epsilon)$-spanning set }
\right\}
$$
and 
$$
Z_1(n, \epsilon) = \sup\left\{
\sum_{y \in \Sigma } \exp(G^n(y))
\hbox{ : } \Sigma  \hbox{ is an $(n, \epsilon)$-{\color{black}separated} set }
\right\}.
$$
where $G^n(x) = \sum_{j=0}^{n-1}G(f^jx)$.
Then 
$$
P(G) = \lim_{\epsilon \to 0}
\lim_{n \to +\infty} \frac{1}{n}\log Z_0(n, \epsilon)
 = \lim_{\epsilon \to 0}
\lim_{n \to +\infty} \frac{1}{n}\log Z_1(n, \epsilon).
$$
\end{lemma}
(See \cite{walters}, chapter 9,  and also  \cite{kh}, Propositions 20.2.7 and 20.3.2)


\begin{figure}[h!]
\centerline{
\begin{tikzpicture}
\draw[fill] (4.05, 2) circle (1.2pt);
\draw[black, ultra thick] (4.05,1.5)--(4.05,2.5);
\node at (5.7,2) {$B_{d_u}(x_i,\epsilon)$};
\draw [domain=-2:-0, samples=50] plot (\x, {2+2*cos(2*pi*\x r)});
\draw (-5,1) .. controls (-4.9,2) .. (-5,3);
\draw [domain= 2:5, samples=150] plot (\x, {2+2*cos(10*pi*\x r)});
\node at (-5.8,1) {$W_\delta^u(x)$};
\node at (-2.7,1) {$f^{m}W_\delta^u(x)$};
\node at (6.3,1) {$f^{n+m}W_\delta^u(x)$};
\draw[fill] (-0.25, 2) circle (1.2pt);
\node at (-0.47, 2){$y_i$};
\draw[fill] (-0.17, 3) circle (1.2pt);
\node at (-0.45, 3){$y$};
\draw[fill] (0.3, 3) circle (1.2pt);
\node at (0.55, 3){$z$};
 \end{tikzpicture}
}
\caption{The  push forward  $f^m W_\delta^u(x)$ is $\epsilon$-dense.}
\end{figure}
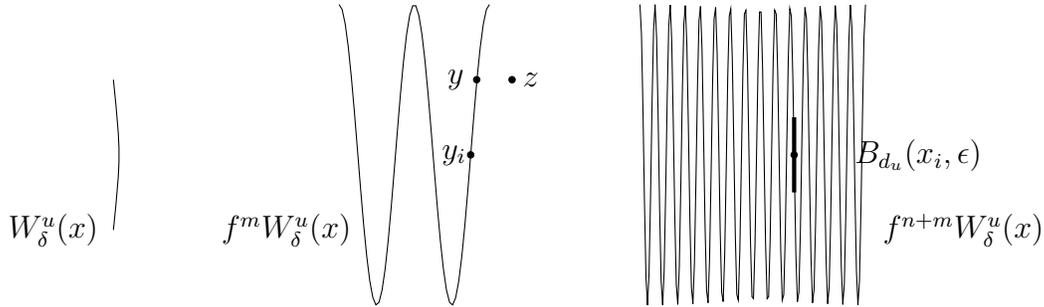

\medskip
To get a lower   bound on the growth rate in Proposition \ref{pressure}
 given $\epsilon > 0$ and 
  $n \geq 1$  we want to 
construct an $(n, 2\epsilon)$-spanning set.  
We begin by choosing a covering of $ f^{n+m} W_\delta^u(x)$
by 
 $\epsilon$-balls 
$$
B_{d_u}(x_i,  \epsilon) \hbox{ : } i=1, \cdots, N:=N(n+m, \epsilon)$$
contained within the unstable manifold with respect to  the induced metric $d_u$
and let  $A_{\epsilon} : = f^{n+m} W_\delta^u(x) \setminus \bigcup_{y \in \partial f^{n+m} W_\delta^u(x)} B_{d_u}(y, \epsilon/2).$
We can choose a maximal set 
$S= \{x_1, \cdots, x_{N(n+m,\epsilon)}\}$ with the additional properties that 
$d_u(x_i,x_j) > \epsilon/2$ for $i \neq j$ and $x_i \in A_\epsilon$.
By our choice of $S$ we have that 
$$
A_{\epsilon} \subset \bigcup_{i=1}^{N(n+m,\epsilon)} B_{d_u}(x_i,  \epsilon/2).
$$
By the triangle inequality we have that 
$$
f^{n+m} W_\delta^u(x) \subset \bigcup_{i=1}^{N(n+m,\epsilon)} B_{d_u}(x_i,  \epsilon).
$$
Since 
$B_{d_u}(x_i, \epsilon/4)\cap B_{d_u}(x_j, \epsilon/4) = \emptyset$ for $i \neq j$
we have that the disjoint union satisfies 
$$
\bigcup_{i=1}^{N(n+m,\epsilon)} B_{d_u}(x_i,  \epsilon/4) \subset f^{n+m} W_\delta^u(x).
$$




We can assume without loss of generality that 
$$f^n: f^{m}W^u_\delta(x)
\to f^{n+m}W^u_\delta(x)$$
 locally expands distance along the unstable 
 manifold (otherwise we can achieve this by a suitable choice of Riemannian metric, e.g.,  the Mather metric).
 In particular, this means that  the primages $y_i := f^{-n} x_i \in  f^m(W_\delta^u(x))$
($i=1, \cdots, N$)
 form an 
$(n, 2 \epsilon)$-spanning set.
By Lemma \ref{cover0},  for  any point  $z\in \textcolor{black}{X}$  we can choose a point
$y \in f^m(W^u_\delta(x))$ with $z \in W^s_\epsilon(y)$
and observe that $d(f^jz, f^jy) < \epsilon$ for $0\leq j \leq n$.
We can then choose a $y_i$ such that $\textcolor{black}{d_n(y,y_i)} < \epsilon$ since $f^n$ is locally expanding along unstable manifolds.
In particular,  by the triangle inequality
$$d(f^j z, f^jy_i)  \leq d(f^j z, f^jy)+d(f^j y, f^jy_i) \leq 2\epsilon$$
 for $0\leq j \leq n$.

Since $G$ is continuous  we have the following bound.

\begin{lemma}\label{bound}
For all  $\tau> 0$ there exists $\epsilon > 0$  sufficiently small such that 
 for all $n \geq 1$ and 
points $y_i, z\in X$ satisfying  $d(f^jy_i, f^jz) \leq \epsilon$ for $0 \leq j \leq n-1$
we have $|G^n(y) - G^n(z)| \leq n \tau$.
\end{lemma}

It remains to relate $Z_0(n, 2\epsilon)$ to an integral over 
$f^{n+m}W_\delta^u(x)$.  
Because of the properties of our choice of 
$\epsilon$-cover for $f^{n+m}W_\delta^u(x)$
we  have that for all $n \geq 1$
$$
\begin{aligned}
Z_0(n,2\epsilon) &\leq \sum_{i=1}^N \exp(G^n(y_i)) 
\cr
&\leq
\sum_{i=1}^N 
\frac{1}{\lambda(B_{d_u}(x_i,\epsilon/4))}
\int_{B_{d_u}(x_i,\epsilon/4)} \exp(G^n(f^{-n}x_i)) d\lambda(\textcolor{black}{z})\cr
&\leq 
\frac{1}{M} e^{n \tau}  \int_{f^{n+m}W_\delta^u(x)}  e^{G^n(f^{-n}\textcolor{black}{z})}d\lambda(\textcolor{black}{z})
\cr
\end{aligned}
$$
where $M = M(\epsilon) = \inf_{z} \lambda(B_{d_u}(z,\epsilon/4))>0$.
Finally,    we can bound 
$$
\int_{f^{n+m}W_\delta^u(x)}  e^{G^n(f^{-n}\textcolor{black}{z})}d\lambda(\textcolor{black}{z})
\leq 
e^{m \|G\|_\infty}
\int_{f^{n+m}W_\delta^u(x)}  e^{G^{n+m}(f^{-(n+m)}\textcolor{black}{z})}d\lambda(\textcolor{black}{z}).
$$
Comparing these inequalities we see that 
$$
\lim_{\epsilon \to  0}
\lim_{n \to +\infty}
\frac{1}{n} \log 
Z_0(n, 2\epsilon) \leq 
\lim_{n \to +\infty}\frac{1}{n} \log \int_{f^{n}W_\delta^u(x)}  e^{G^n(f^{-n}\textcolor{black}{z})}
d\lambda(\textcolor{black}{z}) +\tau.
$$
Since $\tau>0$ can be chosen arbitrarily small the lower bound follows.

To get an upper   bound on the growth rate  in Proposition \ref{pressure}, 
given $\epsilon > 0$ and 
  $n \geq 1$ then we want to
create an $(n,\kappa \epsilon)$-separated  set.  
To this end, we  can  choose a maximal number of points $x_i \in \textcolor{black}{f^n}W_\delta^u(x)$ $(i = 1, \cdots, N= N(n, \epsilon))$
so that $d_u(x_i, x_j) > \epsilon$ whenever $i \neq j$. 
 We can again assume without loss of generality that 
 $f^n: W^u_\delta(x) \to f^nW^u_\delta(x)$ is locally  distance expanding and thus, in particular, the points
 $y_i = f^{-n}x_i$ 
 ($i=1, \cdots, N = N(n, \epsilon)$)
 form an $(n, \kappa \epsilon)$-separated set, for some $\kappa > 0$ independent of $n$ and $\epsilon$. 
 
The balls 
$B_{d_u}(x_i,  \epsilon)$ ($i=1, \cdots, N=N(n, \epsilon)$) form a cover for $f^nW^u_\delta(x)$, since otherwise we could choose an extra point  $z \in f^nW^u_\delta(x)$  with   $\inf_i\{d(z, x_i)\} \geq \epsilon$ 
contradicting the maximality of the previous family.
We  can therefore use Lemma 
 \ref{bound} to
 bound
 $$
\begin{aligned}
Z_1(n,\textcolor{black}{\kappa} \epsilon) &\geq \sum_{i=1}^N
\frac{e^{-n\tau}}{\lambda(B_{d_u}(x_i, \epsilon))}
 \int_{B_{d_u}(x_i,\epsilon)}  \exp({G^n(f^{-n}z))}d\lambda(z)\cr
&\geq 
\frac{e^{-n\tau }}{L}\int_{f^{n}W_\delta^u(x)}  e^{G^n(f^{-n}\textcolor{black}{z})}d\lambda(\textcolor{black}{z}).  
\end{aligned}
$$
 where  $L = L(\epsilon)= \sup_z  \lambda(B_{d_u}(z, \epsilon)) > 0$.
In particular,  we see that 
$$
\lim_{\epsilon \to 0}
\lim_{n \to +\infty}
\frac{1}{n} \log 
Z_1(n, \kappa \epsilon) \geq 
\lim_{n \to +\infty}\frac{1}{n} \log \int_{f^{n}W_\delta^u(x)}  e^{G^n(f^{-n}\textcolor{black}{z})}d\lambda(\textcolor{black}{z}) - \tau.
$$
Since $\tau > 0$ is arbitrary 
this completes the proof.

 \end{proof}

\begin{remark}
Under the  stronger assumption that $G$ is H\"older continuous then in the above argument we can replace $n\tau$ by a constant $C>0$ independent of $n$ by using a telescoping argument.  This leads to some simplifications.
\end{remark}
 
 We can consider 
 Proposition \ref{pressure} in the special case of 
 the two distinguished examples of the SRB measure and the measure of maximal entropy.

\begin{example}
In the particular case that  $G= \Phi$ then we see that 
$$  \int_{W_\delta^{u}(x)}
\underbrace{
 \exp\left( \sum_{k=0}^{n-1} (G-\Phi)(f^kz) \right) }_{=1}
 d\textcolor{black}{\lambda}(z)=
 \textcolor{black}{\lambda}(W_\delta^{u}(x))
 $$
 which is independent of $n$ and thus  we recover from Proposition \ref{pressure}
 that 
 $P(\Phi)=0$ (see \cite{BR}).
\end{example}

\begin{example}
In particular, when $G= 0$ then we see that the topological entropy can be written in the form 
$$  
 h(f) = \lim_{n \to +\infty}
 \frac{1}{n} \log 
 \left(
  \lambda_{f^nW_\delta^u}(f^nW_\delta^{u}(x))
  \right).
 $$
\end{example}

There have been previous results related to growth rates of the lengths of curves.
For general $C^\infty$  diffeomorphisms, for example,  Newhouse has an expression for the entropy in terms of the growth of pieces of unstable manifolds \cite{newhouse}.  

\section{Proof of Theorem \ref{thm:gibbsdiffeo}}
By Alaoglu's theorem on the weak star compactness of the 
space of $f$-invariant
probability measures
we can find an $f$-invariant  probability measure, which we denote by   $\mu$,  and 
a  subsequence $n_k$ such that 
the measures  (1.\textcolor{black}{3}) have a weak star convergent subsequence
$\lim_{k \rightarrow \infty} \mu_{n_k} = \mu$.
Moreover, for  any continuous $F: X \to \mathbb R$ we can compare 
$$
\begin{aligned}
|\int F d\mu_n - \int F\circ  f d\mu_n| &= 
|\frac{1}{n}\sum_{k=0}^{n-1} \int F\circ f^k d\lambda_n
- \frac{1}{n}\sum_{k=0}^{n-1} \int F\circ f^{k+1} d\lambda_{n}|\cr
&\leq \frac{2\|F\|_\infty}{n} \to 0 \hbox{ as } n \to +\infty\cr
\end{aligned}
$$
and, in particular,  one easily sees that $\mu$ is $f$-invariant.

For convenience we denote
$$
\begin{aligned}
Z_n^G 
&= \int_{W_\delta^u(x)} \exp \left( 
\sum_{k=0}^{n-1} (G-\Phi)(f^ky)
\right)d\lambda_{W_\delta^u(x)}(y)\cr
&
=  \int_{f^nW_\delta^u(x)} \exp \left( 
\sum_{k=1}^{n} G(f^{-k}y)
\right)d\lambda_{f^nW_\delta^u(x)}(y)
\cr
\end{aligned}
$$
using the change of variables on $f^n: W_\delta^u(x) \to f^nW_\delta^u(x)$.
We want to show that
$\mu$ is the Gibbs measure for $G$.
As $G$ is H\"older continuous then  by uniqueness of the Gibbs measure $\mu_G$ we would  have that $\mu_n \to \mu_G$
as  $n \to +\infty$.

\begin{definition}
Given a finite partition $\mathcal P = \{P_i\}_{i=1}^N$ we say that it has size $\epsilon > 0$ if $\sup_{i}\left\{ \hbox{\rm diam}(P_i)  \right\}< \epsilon$.
\end{definition}




\textcolor{black}{
Since unstable manifolds are $C^1$, there exists a constant $C > 0$ such that for all $y,z \in f^{n-1} W_\delta^u(x)$ and for $\epsilon > 0$ sufficiently small, if $d_u(y,z) > \epsilon$ then there exists an $0 \leq i \leq n-1$ such that $d(f^{-i} y, f^{-i} z) > C\epsilon$. In particular, $d_{n}(f^{-(n-1)}y, f^{-(n-1)}z) > C\epsilon$.}
\textcolor{black}{By Lemma 4.9}, for  any $\tau>0$ we can choose a partition $\mathcal P$  of size 
\textcolor{black}{$C \epsilon$} 
  such that  for all  $x,y \textcolor{black}{ \in A \in} \vee_{i=0}^{n-1}
f^{-i}\mathcal P$ we have that
$$\left|\sum_{k=0}^{n-1} G(f^{k}x) - \sum_{k=0}^{n-1} G(f^{k}y)\right| \leq n \tau. \eqno(5.1)$$

{\color{black} Under the assumption that $G$ is H\"older continuous the upper bound in (5.1) can be replaced by $\tau$, say. However, if we only assume that $G$ is continuous then the above bound suffices and then the proof also applies to give Theorem \ref{thm:gibbsdiffeo-1}.
}



Proceeding with the proof of  Theorem \ref{thm:gibbsdiffeo}, 
for each $A \in  \bigvee_{h=0}^{n-1} f^{-h} \mathcal P$ we can fix \textcolor{black}{an} $x_A \in A$.
By definition of $\lambda_n$, 
 $$
\begin{aligned} 
\int_{W^u_\delta(x)} G(x) d \lambda_{n}(x)& = 
  \frac{1}{Z_n^G}  \int_{f^nW^u_\delta(x)} 
\exp\left({\sum_{k=1}^{n}G(f^{-k}y)}\right)G(f^{-n}y)  d\lambda_{f^nW^u_\delta(x)}(y)
\end{aligned}
$$
and so we can write for each $0 \leq j \leq n-1$:
$$
\begin{aligned}
\int_{f^jW^u_\delta(x)} G d (\textcolor{black}{f^j_*}\lambda_{n}) &= \int_{W_\delta^u(x)} G( f^{j} y) d \lambda_{n}(y) \cr
  &= \frac{1}{Z_n^G}  \int_{f^n(W^u_\delta(x))} 
\exp\left({\sum_{k=1}^{n} G(f^{-k}y)}\right)G( f^{-(n-j)}y)  d\lambda_{f^nW^u_\delta}(y).
\end{aligned}
$$
Hence,
using  the definition of $\mu_n$  in (1.\textcolor{black}{3}) we can write
$$
\begin{aligned}
&\int G(y) d\mu_{n}(y) \cr
 & = \frac{1}{nZ_n^G}
 \int_{f^nW^u_\delta(x)} 
\exp\left({\sum_{k=1}^{n}G(f^{-k}y)}\right)
\left( \sum_{j=0}^{n-1}G( f^{-(n-j)}y) \right) d\lambda_{f^nW^u_\delta}(y)\\
			& = \frac{1}{nZ_n^G} \sum_{A \in \bigvee_{h=0}^{n-1} f^{-h}\mathcal P} \int_{f^n(A \cap W^u_\delta(x))} 
\exp\left({\sum_{\textcolor{black}{k=1}}^{\textcolor{black}{n}}G(f^{-k}y)}\right)
\left( \sum_{j=0}^{n-1}G( f^{-(n-j)}y) \right) d\lambda_{f^nW^u_\delta}(y)\\			
			& \geq \frac{1}{nZ_n^G} \sum_{A \in \bigvee_{h=0}^{n-1} f^{-h}\mathcal P} 
			\left( \sum_{j=0}^{n-1}G( f^{j}x_A) - n\tau \right)    \int_{f^n(A \cap W^u_\delta(x))} \exp\left({\sum_{k=1}^n G(f^{-k}y)}\right)  d\lambda_{f^nW^u_\delta}(y)
			 \label{pressureeq1}
\end{aligned}
\eqno(5.2)
$$
the last inequality coming from equation (5.1). 
We can write for a set $A \in \bigvee_{h=0}^{n-1} f^{-h}\mathcal P$:
$$
\lambda_n(A) = \int_{f^nW^u_\delta(x)} 
\exp\left({\sum_{k=1}^{n}G(f^{-k}y)}\right)
\chi_{f^nA}(y) d\lambda_{f^nW^u_\delta(x)} (y), 
$$
then 
writing $\lambda_{f^n (A\cap W^u_\delta(x))}$ for the restriction 
of $\lambda_{f^n (W^u_\delta(x))}$  to $f^n A$
we have that 
$$
\begin{aligned}
\log \lambda_n(A)
					& =
		 \log  \int_{f^nW^u_\delta(x)} 
\exp\left({\sum_{\textcolor{black}{k=1}}^{\textcolor{black}{n}}G(f^{-k}y)}\right)
 d\lambda_{f^n (A\cap W^u_\delta(x))} 		\cr	
&\leq \sum_{k=0}^{n-1}G(f^{k}x_A) + n\tau  + \log \lambda_{f^n W^u_\delta(x))}(f^n(A))
\cr
&\leq \sum_{k=0}^{n-1}G(f^{k}x_A) + n \tau 
\end{aligned}
\eqno(5.3)
$$

\textcolor{black}{
The last inequality requires the following justification.
}

\textcolor{black}{
First note that as $A \in \bigvee_{h=0}^{n-1} f^{-h}\mathcal P$ and since $\mathcal{P}$ is of size $C \epsilon$ then by definition, there is a $z \in X$ such that $A \in B_n(z, C\epsilon)$. We will now show, for any $\epsilon' >0$, there is a $p \in f^nW_\delta^u(x)$ and $\epsilon > 0$ small enough such that $f^nA \cap f^nW_\delta^u(x) \subset B_{d_u}(p,\epsilon')$ and in particular, the intersection is contained within one small piece of local unstable manifold, $f^nA \cap f^nW_\delta^u(x) \subset W_{\epsilon'}^u(p)$.
}

\textcolor{black}{
By continuity, for $\epsilon > 0$ small enough, if $f^{n-1}A \cap f^{n-1}W_\delta^u(x) \subset B_{d_u}(f^{-1}p,\epsilon)$ then $f^nA \cap f^nW_\delta^u(x) \subset B_{d_u}(p,\epsilon')$. Therefore, we only need to check $f^{n-1}A \cap f^{n-1}W_\delta^u(x) \subset B_{d_u}(f^{-1}p,\epsilon)$ for some $f^{-1}p \in f^{n-1}W_\delta^u(x)$. If this were not the case then there exist $z_1, z_2 \in f^{n-1}A \cap f^{n-1}W_\delta^u(x)$ such that $B_{d_u}(z_1, \frac{\epsilon}{2}) \cap B_{d_u}(z_2, \frac{\epsilon}{2}) = \emptyset$ (since $A \in \bigvee_{h=0}^{n-1} f^{-h}\mathcal P$ and $\mathcal{P}$ is of size $C\epsilon$). Then $d_u(z_1,z_2) > \epsilon$ and so $d_n(f^{-(n-1)}z_1, f^{-(n-1)}z_2) > C\epsilon$ (by the choice of $C$ above equation (5.1)), but this contradicts that $A \in B_n(z, C\epsilon)$ since $f^{-(n-1)}z_1, f^{-(n-1)}z_2 \in A$. 
 }

\textcolor{black}{
Therefore, as $f^nA \cap f^nW_\delta^u(x) \subset B_{d_u}(p,\epsilon')$, we can take $\epsilon'$ sufficiently small (by taking $\epsilon$ small enough) such that $\lambda_{f^n W^u_\delta(x))}(f^n(A)) < 1$ and so $\log \lambda_{f^n W^u_\delta(x))}(f^n(A))$ is negative as required.
}




Letting $K_{n,A} =   \int_{f^nW^u_\delta(x)} 
\exp\left({\sum_{k=1}^{n}G(f^{-k}y)}\right)
 d\lambda_{f^n (A\cap W^u_\delta(x))} $
 we can 
 consider the entropy
$$
\begin{aligned}
H_{\lambda_{n}} \bigg(\bigvee_{h=0}^{n-1} f^{-h}\mathcal{P} \bigg) &
 = - \sum_{A \in \bigvee_{h=0}^{n-1} f^{-h}\mathcal{P}} \lambda_{n}(A) \log \lambda_{n}(A), \\
			& = - \sum_{A \in \bigvee_{h=0}^{n-1} f^{-h}\mathcal P} \frac{K_{n,A}}{Z_n^G} \log \frac{K_{n,A}}{Z_n^G}, \\
			& = \log Z_n^G - \sum_{A \in \bigvee_{h=0}^{n-1} f^{-h}\mathcal P} \frac{K_{n,A}}{Z_n^G} \log{K_{n,A}}
\end{aligned}
\eqno(5.4)
$$
where the last equality uses
$\sum_{A \in \bigvee_{h=0}^{n-1}f^{-h}\mathcal P} K_{n,A} =  Z_n^G.$
Therefore,  comparing  (5.3) and  (5.4) gives
$$
H_{\lambda_{n}} \bigg(\bigvee_{h=0}^{n-1} f^{-h}\mathcal P \bigg) \geq \log Z_n^G  - \sum_{A \in \bigvee_{h=0}^{n-1} f^{-h}\mathcal P} \frac{K_{n,A}}{Z_n^G}
\left(\sum_{k=0}^{n-1}G(f^{k}x_A)+  n \tau\right).
\eqno(5.5)
$$
By  (5.2) we can also bound
$$
n \int_X G d \mu_{n} \geq 
  \frac{1}{Z_n^G} \sum_{A \in \bigvee_{h=1}^{n} f^{-h}\mathcal P} \left(\sum_{k=1}^nG(f^{k}x_A) -  n \tau\right) K_{n,A}.
  \eqno(5.6)
$$
{\color{black} Summing} (5.5) and (5.6) we can write
$$
\begin{aligned}
&H_{\lambda_{n}} \bigg(\bigvee_{h=0}^{n-1} f^{-h}\mathcal P \bigg) + n \int_X G d \mu_{n} 
\cr& \geq \log Z_n^G - \sum_{A \in \bigvee_{h=0}^{n-1} f^{-h}\mathcal{P}} \frac{K_{n,A}}{Z_n^G}\left(\sum_{k=0}^{n-1}G(f^{k}x_A) 
+ n \tau\right)\cr
&\qquad\qquad\qquad
+ \frac{1}{Z_n^G} \sum_{A \in \bigvee_{h=0}^{n-1} f^{-h}\mathcal P} \left(\sum_{k=0}^{n-1}G(f^{k}x_A)-  n \tau\right) K_{n,A}, \\
			& \geq \log Z_n^G  - 2 n \tau.
\end{aligned}
\eqno(5.7)
$$
We can use (5.7) and Lemma 4.5 to write,
$$
\begin{aligned}
q \log Z_n^G - q n \int_X G d \mu_{n}-2q n \tau & \textcolor{black}{\leq}
q H_{\lambda_{n}} \bigg(\bigvee_{h=0}^{n-1} f^{-h}\mathcal{P} \bigg) , \\ 
				&	 \leq n H_{\mu_{n}} \bigg(\bigvee_{i=0}^{q-1} f^{-i}\mathcal{P} \bigg) + 2q^2
				{\color{black} \log} \hbox{\rm Card}(\mathcal{P}), 
\end{aligned}
$$
{\color{black} then dividing by $q$ and $n$, and rearranging the terms gives}
$$
\begin{aligned}				
\frac{\log Z_n^G}{n} - \frac{2n\tau}{n} - \frac{2q {\color{black} \log}|\mathcal{P}|}{n} & \leq \frac{H_{\mu_{n}} \bigg(\bigvee_{i=0}^{q-1} f^{-i}\mathcal{P} \bigg)}{q} + \int_X G d \mu_{n}.
\end{aligned}
$$
Letting $n_k \to +\infty$ 
\begin{align*}
P(G) & = \lim_{k \rightarrow \infty} \frac{\log Z_{n_k}^G}{n_k}\cr
& \leq \lim_{k \rightarrow \infty} \bigg(\frac{H_{\mu_{n_k}}\bigg(\bigvee_{i=0}^{q-1} f^{-i}\mathcal{P} \bigg)}{q} + \int_X G d \mu_{n_k}\bigg) + 2 \tau \\
			& =  \frac{H_{\mu}\bigg(\bigvee_{i=0}^{q-1} f^{-i}\mathcal{P} \bigg)}{q} + \int_X G d \mu + 2\tau,
\end{align*}
where we assume without loss of generality that the boundaries of the partition have zero measure.
Letting  $q \rightarrow \infty$,
$$
\begin{aligned}
P(G) \leq \textcolor{black}{h_{\mu}(\mathcal P)} + \int_X G d \mu + \textcolor{black}{2}\tau.
\end{aligned}\eqno(5.8)
$$
Finally, we recall that  $\tau$ is arbitrary.
Therefore, since $\mu$ is a $f$-invariant probability measure 
we see from the variational principle (1.1) that the inequalities in (5.8) are actually equalities 
(since \textcolor{black}{$h_{\mu}(\mathcal P)  \leq  h({\mu})$})
and therefore  we conclude  that the
measure  $\mu$ is the  Gibbs measure for $G$, $\mu_G$.

\section{Proof of Theorem \ref{thm:gibbsflow}}


The proof of  Theorem \ref{thm:gibbsflow} for Gibbs measures for  flows is  completely analogous to that for diffeomorphisms in Theorem \ref{thm:gibbsdiffeo}.

\subsection{Hyperbolic attracting  flows}
We begin by recalling the definition of a hyperbolic attracting  flow.
Let $\phi_t: M \to M$ 
 ($t \in \mathbb R$)
be a $C^{1+\alpha}$ flow  on a compact \textcolor{black}{Riemannian} manifold, and let $X \subset M$ be a closed $\phi$-invariant set.

\begin{definition}
The flow 
 $\phi_t:X \to X$ is called a  {\it \textcolor{black}{mixing} hyperbolic attractor} if: 
\begin{enumerate}
\item 
there exists a continuous splitting $T_X M =E^0 \oplus E^s\oplus E^u$ 
where $E^0$ is a one dimensional subbundle tangent to the flow orbits 
and  there exist $C > 0$ and $0< \lambda < 1$ such that 
$$\|D\phi_t|E^s\| \leq C \lambda^t \hbox{ and }\|D\phi^{-t}|E^u\| \leq C \lambda^t$$ for $t \geq 0$;
\item there exists an open set $X \subset U \subset M$ such that 
$X = \cap_{t=0}^\infty \phi_tU$;
\item
$\phi_t: X \to X$ is topologically mixing; and
\item
the periodic orbits for 
$\phi_t: X \to  X$ are dense in $X$.
\end{enumerate}

In the particular case that the entire manifold is hyperbolic, i.e., $X = M$ then we call the flow Anosov. 
\end{definition}




  Given any H\"older continuous function $G : X \to \mathbb R$  we want to consider an analogous construction of a Gibbs measure.    


 Let us first define a H\"older continuous function 
  $\Phi: \textcolor{black}{X} \to \mathbb R$ defined by 
  $$
  \Phi(y) =  - \lim_{t\to 0} \frac{1}{t} \log |\textcolor{black}{\hbox{\rm det} (D\phi_t | E_y^u)}|.
  $$

\begin{example}
\begin{enumerate}
\item
In the particular case that 
$G=0$ we have that the limit is 
the measure of maximal entropy.
\item
In the particular case that 
$G=  \Phi(x)$ we have that the limit is 
$\mu_{SRB}$, the SRB measure.
\end{enumerate}
\end{example}



The proof of Theorem 
\ref{thm:gibbsflow}
follows the same lines as Theorem  \ref{thm:gibbsdiffeo}.
We will therefore only need to explain the main ideas.

\subsection{Basic properties} 
We next recall the definition of the stable and unstable manifolds for the flow. 
  
 \begin{definition}
We associate to each $x \in X$ the  {\it stable manifold} defined by
$$
W^{s}(x) = \{y \in M \hbox{ : } d(f\phi_tx, \phi_ty) \to 0 \hbox{ as } t \to +\infty\}
$$
and  {\it unstable manifold}  defined by
$$
W^{u}(x) = \{y \in M \hbox{ : } d(\phi_{-t}x, \phi_{-t}y) \to 0  \hbox{ as } t \to +\infty\}.
$$
For sufficiently small $\delta > 0$ we define local versions 
$W_\delta^{s}(x) \subset W^{s}(x)$ and  $W_\delta^{u}(x) \subset W^{u}(x)$
defined by 
$$
W_\delta^{s}(x) =\{y \in M \hbox{ : } d(\phi_tx, \phi_ty) \to 0 \hbox{ as } t \to +\infty 
\hbox{ and }d(\phi_tx, \phi_ty)  \leq \delta, \forall t \geq 0\}
$$
and 
$$
W_\delta^{u}(x) =  \{y \in M \hbox{ : } d(\phi_{-t}x, \phi_{-t}y) \to 0  \hbox{ as } t \to +\infty
\hbox{ and } d(\phi_{-t}x, \phi_{-t}y) \leq \delta, \forall t \geq 0\}.
$$
\end{definition}

The basic properties of the (local) stable and unstable manifolds is the following
analogue of Lemma \ref{stab}

\begin{lemma}
For each $x \in \Lambda$ and $\delta > 0$ sufficiently small:
\begin{enumerate}
\item
the sets
$W^{s}(x)$ and $W^{u}(x)$
are $C^{1}$ immersed submanifolds of dimension 
$\dim (E^s)$ and $\dim (E^u)$, respectively, with $T_xW^s = E_x^s$
and $T_xW^u = E_x^u$.
\item
 the sets
$W_\delta^s(x)$ and $W_\delta^{u}(x)$
are $C^{1}$ embedded  disks  of dimension 
$\dim (E^s)$ and $\dim (E^u)$, respectively, with $T_xW_\delta^s(x) = E_x^s$
and $T_xW^u_\delta(x) = E_x^u$.
\end{enumerate}
\end{lemma}

\subsection{Pressure and growth rates for hyperbolic attracting  flows}

 We will need the following useful characterization of the pressure from the flow, which is analogous to the corresponding result for diffeomorphisms (i.e., Proposition \ref{pressure}).

\begin{proposition}\label{pressure1}
Let $\phi_t: X \to X$ be a \textcolor{black}{mixing}  attracting hyperbolic  flow. 
For any  continuous function $G: X  \to \mathbb R$
$$
P(G)
= \lim_{t \to +\infty} \frac{1}{t}\log  \int_{W_\delta^{u}(x)}
 \exp\left( \int_0^t  (G-\Phi)(\phi_vx)  dv\right) 
 d\lambda_{W_\delta^u(x)}(x).$$
\end{proposition}
Again, we can  use the change of variables to rewrite this as
$$
P(G)
= \lim_{t \to +\infty} \frac{1}{t}\log  \int_{\phi_tW_\delta^{u}(x)}
 \exp\left( \int_0^t  G(\phi_{-v}x)  dv\right) 
 d\lambda_{\phi_tW_\delta^u(x)}(x).$$

\begin{example}
In particular, when $G= \Phi$ then we see that 
$$  \int_{W_\delta^{u}(x)}
\underbrace{
 \exp\left( \int_0^t  (G-\Phi)(\phi_vx) dv \right) }_{=1}
 d\lambda(x)=
 \lambda(W_\delta^{u}(x))
 $$
 which is independent of $t$ and thus  we recover from Proposition 
 6.5
that  $P(\Phi)=0$.
\end{example}

The   proof of Proposition \ref{pressure1} is  analogous to the corresponding result for hyperbolic  diffeomorphisms \textcolor{black}{so we will only provide details where the proof differs}.

 \begin{proof}
 We begin with the following analogue of Lemma \ref{cover0}.

\begin{lemma}\label{cover}
Let $\phi_t: X \to X$ be a mixing hyperbolic  flow.
For any $\epsilon > 0$ there exists $T_0 > 0$ such that $\phi_{T_0}W^u_\delta(x)$ 
is $2 \epsilon$-dense.  In particular, we can assume that
 $$\textcolor{black}{X} = 
 \phi_{[-\epsilon, \epsilon]}
\left( \cup_{y \in \phi_{T_0} W^u_\delta(x)} W_\epsilon^s(y)\right).
 $$
\end{lemma}

\begin{proof}
This is a consequence of the minimality of the unstable foliation.  
\end{proof}


The pressure $P(G)$ for flows can   be written in terms of the growth rates of sets of spanning sets and separated sets. 
Recall that a 
 $(T,\epsilon)$-spanning set $S \subset X$ 
is such that $ \cup_{x\in S} B(x,T,\epsilon)$ covers $X$,  where $B(x,T,\epsilon) := \cap_{0 \leq t \leq T} \phi_{-t} B(\phi_tx, \epsilon)$  
is a   Bowen ball.
On the other hand  a $(T,\epsilon)$-separated  set $\Sigma \subset X$ 
is  one such that 
$\max_{0 \leq t \leq T} d(\phi_tx,\phi_ty) > \epsilon$ for $x,y \in \Sigma$
($x\neq y$).
(In particular, the sets 
$B(x,T,\epsilon/2)$  are disjoint
in $X$).

\begin{lemma} 
Given 
$n \geq 1$ and $\epsilon > 0$ we denote 
$$
Z_0(T, \epsilon) = \inf\left\{
\sum_{x \in S} \exp(G^T(x))
\hbox{ : } S \hbox{ is a $(T, \epsilon)$-spanning set }
\right\}
$$
and 
$$
Z_1(T, \epsilon) = \sup\left\{
\sum_{x \in \Sigma } \exp(G^T(x))
\hbox{ : } \Sigma  \hbox{ is a $(T, \epsilon)$-separated set }
\right\}.
$$
where $G^T(x) = \int_0^TG(\phi_tx) dt$.
Then 
$$
P(G) = \lim_{\epsilon \to 0}
\lim_{T \to +\infty} \frac{1}{T}\log Z_0(T, \epsilon)
 = \lim_{\epsilon \to 0}
\lim_{T \to +\infty} \frac{1}{T}\log Z_1(T, \epsilon).
$$
\end{lemma}





To get a lower bound on the growth rates in the statement of the Proposition, 
 given $\epsilon > 0$ and 
  $T > 0$ then we want to
construct  a $(T, \textcolor{black}{3}\epsilon)$-spanning set.  


The construction is analogous with the diffeomorphism case.
We begin by choosing a covering of $ \phi_{T+T_0} W_\delta^u(x)$
by 
 $\epsilon$-balls 
$$
B_{d_u}(x_i,  \epsilon) \hbox{ : } i=1, \cdots, N:=N(n+m, \epsilon)$$
contained within the unstable manifold with respect to  the induced metric $d_u$
and let  $A_{\epsilon} : =  \phi_{T+T_0} W_\delta^u(x) \setminus \bigcup_{y \in \partial\phi_{T+T_0}  W_\delta^u(x)} B_{d_u}(y, \epsilon/2).$
We can choose a maximal set 
$S= \{x_1, \cdots, x_{N(n+m,\epsilon)}\}$ with the additional properties that 
$d_u(x_i,x_j) > \epsilon/2$ for $i \neq j$ and $x_i \in A_\epsilon$.
By our choice of $S$ we have that 
$$
A_{\epsilon} \subset \bigcup_{i=1}^{N(n+m,\epsilon)} B_{d_u}(x_i,  \epsilon/2).
$$
By the triangle inequality we have that 
$$
\phi_{T+T_0}  W_\delta^u(x) \subset \bigcup_{i=1}^{N(n+m,\epsilon)} B_{d_u}(x_i,  \epsilon).
$$
Since 
$B_{d_u}(x_i, \epsilon/4)\cap B_{d_u}(x_j, \epsilon/4) = \emptyset$ for $i \neq j$
we have that the disjoint union satisfies 
$$
\bigcup_{i=1}^{N(n+m,\textcolor{black}{\epsilon})} B_{d_u}(x_i,  \textcolor{black}{\epsilon/4}) \subset\phi_{T+T_0} W_\delta^u(x).
$$


We can assume without loss of generality that 
$$\phi_T: \phi_{T_0} W^u_\delta(x)
\to \phi_{T_0+T} W^u_\delta(x)$$
 locally expands distance along the unstable 
 manifold (otherwise we can achieve this by a suitable choice of Riemannian metric, eg,  the Mather metric).
 
We claim  that  the primages $y_i := \phi_{-T} x_i \in  \phi_{T_0}(W_\delta^u(x))$
($i=1, \cdots, N$)
 form a 
$(T, 3 \epsilon)$-spanning set.
By Lemma \ref{cover} for  any point  $\textcolor{black}{\hat{z}} \in \textcolor{black}{X}$  we can choose a point \textcolor{black}{$z = \phi_{\hat{t}} \hat{z}$ for $\hat{t} \in [-\epsilon, \epsilon]$ and a point} 
$y \in \phi_{T_0}(W^u_\delta(x))$ with $z \in W^s_\epsilon(y)$
and observe that $d(\phi_tz, \phi_ty) < \epsilon$ for $0\leq t \leq T$.
We can then choose a $y_j$ such that \textcolor{black}{$d(\phi_ty,\phi_ty_j) < \epsilon$ for $0\leq t \leq T$} since $\phi_t$ is locally expanding along unstable manifolds.
In particular,  by the triangle inequality 
$$\textcolor{black}{d(\phi_t \hat{z}, \phi_t y_j)  \leq d(\phi_t \hat{z}, \phi_t z) + d(\phi_t z, \phi_ty)+d(\phi_t y, \phi_t y_i) \leq 3 \epsilon}$$
 for $0\leq t \leq T$.

Since $G$ is continuous  we have the following bound.

\begin{lemma}\label{bound1}
For all  $\tau> 0$ there exists $\epsilon > 0$  sufficiently small such that 
 for all $T \geq 1$ and 
points $y,z \in X$ satisfying  $d(\phi_ty,  \phi_tz) \leq \epsilon$ for $0 \leq t \leq T$
we have $|G^T(y) - G^T(z)| \leq T \tau$.
\end{lemma}

It remains to relate $Z_0(T, \textcolor{black}{3}\epsilon)$ to an integral over $\phi_{T+T_0}W_\delta^u(x)$\textcolor{black}{, which is analogous to the diffeomorphism case.}

\textcolor{black}{To get an upper bound on the growth rate in Proposition 6.5 we proceed analogously to the diffeomorphism case to conclude the proof.}


 \end{proof}


\begin{remark}\label{rem:other}
  One can see from the proof that in the statements of 
  Proposition \ref{pressure1} and 
  Theorem \ref{thm:gibbsflow}
  extend to the case that $W^u_\delta(x)$  is replaced by any  embedded submanifold in dimension $\dim E^u$, provided it is not contained in a weak stable manifold.
  \end{remark}

\subsubsection{Growth rates and  geodesic flows on surfaces of variable  negative curvature}

We want to relate the results on hyperbolic attracting  flows to classical results on geodesic flows.

Given a compact manifold $V$ with negative sectional curvatures we can associate the 
unit tangent bundle $M = SV$.  Let $\phi_t: M \to M$ be the geodesic flow.   
Given a point $x \in V$  we can denote by $S_xV$  the fibre above $x$.  
We can consider the image  $\phi_t(S_xV)$ under the geodesic flow for time $t>0$.

We begin with the results on entropy and pressure due to Manning and Ruelle.

\begin{thm}[after Manning] 
We can write
$$
h(\phi) = \lim_{t\to +\infty} \frac{1}{t}\log \lambda (\phi_t(S_xV))
$$
where $\lambda$ denotes the induced volume.
\end{thm}

The original statement of Manning was for the rate of growth of the volume of a ball in the universal cover.  However, because of the hypothesis of negative curvature this corresponds to the rate of growth of the induced  volume of the boundary.

\begin{remark}
It is easy to see that the above result is closely related to 
Remark \ref{rem:other} and Proposition \ref{pressure1}. 
 In particular, if we partition $S_xV$ into a finite number of pieces then we see that the images of these pieces under the geodesic flow can be uniformly approximated by longer pieces of unstable manifold.  This point of view is similar to the classical approach of Marcus to horocycles \cite{marcus}.
\end{remark}

Given a  continuous function $G: M \to \mathbb R$  there is a natural generalization which  takes the following form.
Let $G_{(x,y)} = \int_{[x,y]} G$ be the integral of $G$ along the canonical geodesic segment from $x$ to $y$ on $V$.

\begin{thm}[after Ruelle]\label{thm:ruelle}
We can write
$$
P(G)  
= \lim_{t\to +\infty} \frac{1}{t}
\log \int_{\phi_t(S_xV)) } \exp\left( G_{(x,y)}\right) d \lambda_{\phi_t(S_xV)}(y) 
$$
where $\lambda_{\phi_t(S_xV)}$ denotes the induced volume on $\phi_t(S_xV)$.
\end{thm}

The original statement of Ruelle was a generalization of the result  of Manning 
and formulated on the universal cover $\widetilde M$. The above version can be deduced as before.

\subsection{Entropy and  the end of the  proof of Theorem \ref{thm:gibbsflow}}

The arguments in \S 3 did not require that the transformation is hyperbolic
merely a homeomorphism.  In particular, we can apply the same reasoning where 
{\color{black} we let $f= \phi_{t=1}$ be a time one map of a hyperbolic flow.}
In particular,  we can deduce the following analogue of (5.8).
Let $\mu$ denote a  limit point  of $\mu_T$ ($T>0$).

\begin{proposition}
Let $\phi_{t=1}: X \to X$ be a time one \textcolor{black}{map of a} hyperbolic attracting flow and let $G: X \to \mathbb R$ be a 
 continuous function.  Then any limit point $\mu$ of the measures $\mu_T$ satisfies 
$$
P(G) \leq h_{\mu} (\phi_{t=1}, \mathcal P) + \int_X G d\mu
$$
\end{proposition}

We recall that 
$$
h_{\mu} (\phi_{t=1}, \mathcal P)  \leq h_{\mu} (\phi_{t=1}) =h_{\mu} (\phi)
$$
and deduce from the  variational principle  that $\mu$ is the Gibbs measure for $G$, $\mu_G$.
In particular, if $G$ is H\"older continuous then  by uniqueness of the Gibbs measure we see that $\mu_T \to \mu_G$ as $T \to +\infty$.

  \section{Generalizations and   questions}
 
 We can consider some generalizations of these results
by examining the proof of 
Theorem \ref{thm:gibbsdiffeo} \textcolor{black}{and} it is 
possible to formulate  the  result  in a more general setting.

\begin{proposition}\label{prop:abstract} Let $f: M\to M$ be a $C^2$ diffeomorphism and 
let $W \subset M$ be a $C^2$ embedded disk.  Let $G: M \to \mathbb R$ be a 
continuous function and assume that the  analogue of condition (4.4) holds for $W$.  
If we define $\lambda_n$ as above then any  weak star accumulation point
of  $$\mu_n :=\frac{1}{n}\sum_{k=0}^{n-1} f_*^k \lambda_n
$$
will be an equilibrium state for $G$.
\end{proposition}

If we further assume that 
$G$ is H\"older continuous and 
$f$ satisfies expansion and   specification hypotheses  then 
there will be a unique Gibbs measure \cite{hr}.

We can now describe some future direction where this  theorem could be  applied.

\begin{enumerate}
\item 
In the case of general $C^\infty$ surface diffeomorphisms of positive topological entropy
assume  that there exists at least one strong unstable manifold for which
(4.4) holds.    Then this can be used to construct the Gibbs measure.
In the case of the measure of maximal entropy this hypothesis holds by an  observation of Newhouse and Pignataro \cite{np}.
\item 
For a  partially hyperbolic  diffeomorphism
 $f:X \to X$  the analogue of the  SRB measures are naturally replaced by u-Gibbs measures, as introduced by Pesin and Sinai \cite{ps} (see also \cite{dolgopyat}).  
More generally, if $G: \textcolor{black}{X} \to \mathbb R$ is (H\"older) continuous then there still exists Gibbs measures defined using  (1.1) provided that $f$ is $C^\infty$ or that $f$ is entropy expansive (for example, when $\dim E^0 = 1$).
For partially hyperbolic systems uniqueness is more of an issue and often requires some additional hypotheses.
Interesting examples include Quasi hyperbolic toral automorphisms and frame flows.
\item 
For two sided subshifts of finite type $\sigma: \Sigma_A \to \Sigma_A$ we can consider a copy $\{x\}\times \Sigma_A^+$ of the one sided subshift, equipped with a suitable reference measure,  and  
$\Sigma_A^+$ in place of $W_\delta^u(x)$.
\end{enumerate}

It is a natural question to ask about the speed of convergence of 
$\int f \mu_n$ and $\int f d\mu_T$ for hyperbolic  diffeomorphisms and flows, respectively.   This should be related to the speed of mixing or decay of correlations,


\begin{thebibliography}{111}

\bibitem{anosov}
D. V. Anosov, 
{\it 
Geodesic flows on closed Riemannian manifolds of negative curvature},
Proceedings of the Steklov Institute of Mathematics, 1967, 90, 1-235


\bibitem{bowen-periodic}
R. Bowen,  Periodic points and measures for Axiom A diffeomorphisms. Trans. Amer. Math. Soc. 154 (1971), 377–397.

\bibitem{bowen-SLN}
R. Bowen, 
{\it Equilibrium states and the ergodic theory of Anosov diffeomorphisms},
 Lecture Notes in Mathematics 470, Springer, Berlin, 1975 
 
 \bibitem{BR}
 R. Bowen and D. Ruelle,
 The ergodic theory of Axiom A flows.
Invent. Math. 29 (1975)  181–202.

 \bibitem{brin} M. Brin,
  Ergodic theory of frame flows. Ergodic Theory and Dynamical Systems, II (College Park, MD, 1979/1980) (Progress in Mathematics, 21). Birkhäuser, Boston, MA, 1982, pp. 163–183

\bibitem{CLP}
V. Climenhaga, S. Luzzatto and Y. Pesin,
The geometric approach for constructing Sinai-Ruelle-Bowen measures. J. Stat. Phys. 166 (2017) 467–493. 

\bibitem{CPZ-1}
V. Climenhaga, Y. Pesin and A. Zelerowicz,
Equilibrium states in dynamical systems via geometric measure theory. Bull. Amer. Math. Soc. (N.S.) 56 (2019) 569–610.

\bibitem{CPZ-2}
V. Climenhaga, Y. Pesin and A. Zelerowicz,
 Equilibrium measures for some partially hyperbolic systems. 
 J. Mod. Dyn. 16 (2020), 155–205.
 
\bibitem{jl} J.  De Simoi and C. Liverani, 
 Limit theorems for fast-slow partially hyperbolic systems. Invent. Math. 213 (2018) 811–1016.

\bibitem{dolgopyat}
D. Dolgopyat, Lectures on $u$-Gibbs measures, 
http://www2.math.umd.edu/~dolgop/ugibbs.pdf

\bibitem{hr}
N. Haydn and D. Ruelle,
Equivalence of Gibbs and Equilibrium states for homeomorphisms satisfying expansiveness and specification, Commun. Math. Phys. 148 (1992) 155-167

\bibitem{hp}
 M. Hirsch and C. Pugh,  Stable manifolds and hyperbolic sets. 1970 Global Analysis (Proc. Sympos. Pure Math., Vol. XIV, Berkeley, Calif., 1968) pp. 133–163,  Amer. Math. Soc., Providence, R.I.

\bibitem{kh} A. Katok and  B. Hasselblatt, 
{\it Introduction to the Modern Theory of Dynamical Systems}, C.U.P., Cambridge, 1995.

\bibitem{marcus}
B. Marcus, Ergodic properties of horocycle flows for surfaces of negative curvature,  Ann. of Math. (2) 105 (1977) 81–105.

\bibitem{gibbs2}
E. Mihailescu, Approximations for Gibbs states of arbitrary H\"older potentials on hyperbolic folded sets, Discrete \& Cont. Dynam. Sys., 32 (2012) 961-975.

\bibitem{gibbs3}
E. Mihailescu,   Asymptotic distributions of preimages for endomorphisms, Ergod. Th. \& Dynam. Sys. 31 (2011) 9111-934.

\bibitem{misiurewicz}
M. Misiurewicz, A short proof of the variational principle for a ZN+ action on a compact space. International Conference on Dynamical Systems in Mathematical Physics (Rennes, 1975), pp. 147–157. Ast\'erisque, No. 40, Soc. Math. France, Paris, 1976.

\bibitem{newhouse}
S. Newhouse, Continuity properties of entropy. Ann. of Math. (2) 129 (1989) 215–235

\bibitem{np}
S, Newhouse and T. Pignataro,
On the estimation of topological entropy, 
J. Statist. Phys. 72 (1993) 1331–1351.

\bibitem{parry1}
W. Parry,  Bowen's equidistribution theory and the Dirichlet density theorem. Ergodic Theory Dynam. Systems 4 (1984) 117–134.

\bibitem{parry2}
W. Parry,   Equilibrium states and weighted uniform distribution of closed orbits. Dynamical systems (College Park, MD, 1986–87), 617–625, Lecture Notes in Math., 1342, Springer, Berlin, 1988.

\bibitem{ps}
Y. Pesin and Y. Sinai,
Gibbs measures for partially hyperbolic attractors
Erg. Th. $\&$ Dyn. Sys. 2 (1982) 417–438


\bibitem{ruelle}
D. Ruelle,  A measure associated with axiom-A attractors. Amer. J. Math. 98 (1976)  619–654.

\bibitem{sinai}
Ya. Sinai, Markov partitions and Y-diffeomorphisms, Funct. Anal. and Appl.,
2:1 (1968) 64-89.

\bibitem{gibbs1}
Ya Sinai, Gibbs measures in ergodic theory, Russian Math. Surveys, 27 (1972) 21-69.

\bibitem{smale}
S. Smale, Differentiable dynamical systems. Bull. Amer. Math. Soc. 73 (1967) 747--817.

\bibitem{SV} R. Spatzier and D. Visscher, Equilibrium measures for certain isometric extensions of Anosov systems, Ergodic Theory Dynam. Systems, 38 (2018), 1154-1167

\bibitem{walters}
P. Walters,  {\it Ergodic Theory}, Springer, 1982

\bibitem{williams} R. F. Williams, One-dimensional non-wandering sets. Topology 6 (1967), 473–487.

\bibitem{gibbs4} 
L.-S. Young,  What are SRB measures, and which dynamical systems have them?  
J. Stat. Phys. 108 (2002) 733-754.

\end{thebibliography}
\end{document}